\documentclass{amsart}

\usepackage{amsmath}
\usepackage{amsthm}
\usepackage{amsopn}
\usepackage{amssymb}
\usepackage[all]{xy}
\usepackage{graphicx}
\usepackage{lscape}

\newcommand{\abs}[1]{\lvert #1 \rvert}

\newcommand{\br}[1]{\overline{#1}}
\newcommand{\brr}[1]{\overline{\overline{#1}}}
\newcommand{\td}[1]{\widetilde{#1}}
\newcommand{\tdd}[1]{\widetilde{\widetilde{#1}}}
\newcommand{\ZZ}{\mathbb{Z}}

\newcommand{\WW}{\mathbb{W}}
\newcommand{\FF}{\mathbb{F}}
\newcommand{\GG}{\mathbb{G}}
\newcommand{\MS}{\mathbb{S}}

\theoremstyle{definition}
 \newtheorem{thm}{Theorem}[section]
 \newtheorem{cor}[thm]{Corollary}
 \newtheorem{lem}[thm]{Lemma}
 \newtheorem{prop}[thm]{Proposition}
 
 \newtheorem{ex}[thm]{Example}
 
 \newtheorem{rmk}[thm]{Remark}
\newtheorem*{cor*}{Corollary}
\newtheorem*{defn*}{Definition}
\newtheorem*{ex*}{Example}
\newtheorem*{exs*}{Examples}
\newtheorem*{rmk*}{Remark}
\newtheorem*{claim*}{Claim}
\newtheorem*{conventions}{Conventions}
\numberwithin{equation}{section}
\numberwithin{figure}{section}
\DeclareMathOperator{\ext}{Ext}

\DeclareMathOperator{\im}{Im}

\title[On the existence of $v_2^9$]{On the existence of the self map $v_2^9$ on the Smith-Toda complex
$V(1)$ at the prime $3$}
\author{Mark Behrens}
\address{
Department of Mathematics \\
University of Chicago \\
Chicago, IL 60637, U.S.A.}
\author{Satya Pemmaraju}
\address{
Fixed Income Derivatives \\
UBS Warburg \\
Stamford, CT 06901, U.S.A.}
\date{\today}

\begin{document}

\begin{abstract}
Let $V(1)$ be the Smith-Toda complex at the prime $3$.
We prove that there exists a map $v_2^9: \Sigma^{144}V(1) \rightarrow V(1)$
that is a $K(2)$ equivalence.  This map is used to construct various
$v_2$-periodic infinite families in the $3$-primary stable homotopy groups
of spheres.
\end{abstract}

\maketitle

\tableofcontents

\section{Introduction and statement of results}\label{sec:intro}

Let $V(0)$ denote the mod $3$ Moore spectrum.
Let $V(1)$ be the Smith-Toda complex obtained by taking the cofiber of the
self map $v_1: \Sigma^4 V(0) \rightarrow V(0)$ which 
induces multiplication
by $v_1$ in $K(1)$-homology.  This is an example of a type 2 complex.
The periodicity theorem of Hopkins and Smith \cite{HopkinsSmith} 
states that there exists a $v_2$-self map $v: \Sigma^N V(1) \rightarrow V(1)$
which is a $K(2)$-equivalence.  The purpose of this paper is to provide
a minimal such $v_2$-self map.  The main theorem of this paper is stated
below.
\begin{thm}\label{mainthm}
There exists a self-map 
$$ v_2^9: \Sigma^{144} V(1) \rightarrow V(1) $$
whose effect on $K(2)$ homology is multiplication by $v_2^9$.
\end{thm}

The strategy of proving the theorem is straightforward and computational.
We first prove that the element $v_2^9$ in the Adams spectral sequence
(ASS) for computing $\pi_*(V(1))$ is a permanent cycle.  We then prove that
this map extends over $V(1)$.  We use the ASS instead of the Adams-Novikov
spectral sequence (ANSS) because $v_2^9$ is in Adams filtration $9$,
whereas it has Adams-Novikov filtration $0$.  Therefore there are less
potential targets for a differential supported by $v_2^9$ 
in the ASS than in the ANSS. The ASS $E_2$-term is also easier to compute. 

Our method of showing that $v_2^9$ is a permanent cycle in the ASS is to
consider all of the elements of the ASS which could be targets of 
differentials supported by $v_2^9$ which  would not be
detected in the ASS for $eo_2 \wedge V(1)$.
We then will show all of these potential targets either support non-trivial
differentials, or are killed on earlier pages of the spectral sequence.
This requires knowledge of the ASS of $eo_2 \wedge V(1)$, as well as
the $E_2$ term of the ASS converging to
$\pi_*(V(1))$.  

The spectrum $eo_2$ is a connective cover of the spectrum $EO_2$ discussed
in \cite{GoerssHennMahowald}.  This spectrum should be regarded as a
chromatic level $2$ analog to the spectrum $bo$.  The spectrum $eo_2$ is a
ring spectrum, and thus there is a Hurewicz homomorphism
$$ h:V(1) \rightarrow eo_2 \wedge V(1). $$
For our proof of Theorem~\ref{mainthm}, we need to know what the ASS for
$eo_2 \wedge V(1)$ looks like, and what the effect of the Hurewicz
homomorphism is on Adams $E_2$ terms, and this is accomplished in
Section~\ref{sec:eo_2}.  \emph{The reader who like to avoid a digression on
$eo_2$-theory is invited to skip Section~\ref{sec:eo_2} and simply refer to
Figure~\ref{fig:eo_2^V(1)ASS} and Proposition~\ref{prop:eo_2Hurewicz} for
the relevant information.}
All of the methods in this section derive
from unpublished work of Hopkins, Mahowald, and Miller.  In retrospect, 
the resolution we
use to compute $\ext(eo_2\wedge V(1))$ should be compared to
that of Ravenel \cite[ch. 7]{Ravenel}.  

In Section~\ref{sec:E_2} we compute the $E_2$-term of the ASS for $V(1)$
through the $144$-stem (the degree of $v_2^9$).
We rely on 
Tangora's computer generated tables
\cite{Tangora}
of $H^*(P_*)$ where $P_*$ is the polynomial part of the dual
Steenrod algebra.  The periodic lambda algebra \cite{Gray} 
allows us to compute the $E_2$ term $H^*(A_*//E[\tau_0, \tau_1])$ via a
Bockstein spectral sequence (BSS).  In some instances Christian Nassau's
computer generated $\ext$ tables were of welcome assistance.

Differentials in the ASS are computed by using the Hurewicz image in the
homotopy of the spectrum $eo_2$. In addition, we will prove
a modified `Leibnetz rule' for
differentials in the ASS for V(1).  This product rule is our main tool 
for calculating
difficult differentials in the ASS for $V(1)$.
The product rule is presented in Section~\ref{sec:ProductRule}.  It is a
generalization of a formula for Adams $d_2$'s that was 
communicated to us by Brayton
Gray.
Such technology is essential because $V(1)$
is \emph{not} a ring spectrum, so its ASS is not a spectral sequence of
algebras.  However, the $S$-module structure of $V(1)$ does make the ASS a
spectral sequence of modules over the ASS for computing $\pi_*(S)$, and
this is used occasionally to propagate differentials.  

We make heavy use of
the computation of the $3$-primary stable
stems through the $108$ stem presented in \cite{Ravenel}.  
We use these computations as input for 
the Atiyah-Hirzebruch spectral sequence
(AHSS) to make selective computations of $\pi_*(V(1))$ in certain ranges.
These computations are described in Section~\ref{sec:AHSS}

Section~\ref{sec:mainthmpf1} is devoted to proving that $v_2^9$ is a
permanent cycle in the ASS for $V(1)$, and therefore detects an element of 
$\pi_*(V(1))$.

It now remains to extend $v_2^9$ over $V(1)$.
This is simplified by a certain
splitting of the complex $D(V(1))\wedge V(1)$, where $D(V(1))$ is the
Spanier-Whitehead dual.  This splitting is the subject of
Section~\ref{sec:splitting}.  The splitting is also needed to prove the
product rule.

Using the attaching map structure of one of the
wedge summands, the obstruction to the extension of $v_2^9$ over $V(1)$
is identified as an
element of $\pi_*(V(1))$.  In showing that this obstruction is zero, it is
helpful to know what power $\beta_1^k$ has the property that $\beta_1^k:
\Sigma^{10k} V(1) \rightarrow V(1)$ is null.  In 
Section~\ref{sec:beta_1_order}, 
we show that the map is null precisely when $k \ge
5$.

In Section~\ref{sec:mainthmpf2}, we 
proceed by considering all of the elements in the ASS that might survive to
the obstruction to extending $v_2^9$ over $V(1)$, and show that they all 
either support differentials, or
are the targets of differentials.  Thus the map $v_2^9$ extends to a self map 
of $V(1)$, completing our proof of Theorem~\ref{mainthm}.

We will now indicate the construction of some $v_2$-periodic elements of the
stable stems which arise from the self-map $v_2^9$.
The authors learned of these constructions from Katsumi Shimomura (compare
with \cite{Shimomura:beta_1}).
Theorem~\ref{mainthm} allows us to deduce that certain elements of the
ANSS for the sphere must be permanent cycles.  In
particular, we have the following consequence.

\begin{cor}\label{cor:betas}
The elements $\beta_i$ are permanent cycles in the ANSS for $i \equiv 
0, 1, 2, 5, 6 \pmod 9$.
\end{cor}

\begin{proof}
For dimensional reasons, $v_2$ is a permanent cycle of the ANSS which
detects a map
$$ v_2: S^{16} \rightarrow V(1). $$
By \cite{Oka}, or Remark~\ref{rmk:v_2^5}, the element $v_2^5$ in
the ANSS is a permanent cycle which detects a map
$$ v_2^5: S^{80} \rightarrow V(1). $$
We denote the Spanier-Whitehead dual of $v_2$ by 
$$v_2^*: \Sigma^{16}DV(1) = \Sigma^{10}V(1) \rightarrow S^0.$$
Let $\nu: V(1) \rightarrow S^6$ be projection onto the top cell.  The
elements $\beta_i$ in Corollary~\ref{cor:betas} are constructed by the
following compositions.
\begin{align*}
\beta_{9t} & : S^{144t-6} \hookrightarrow \Sigma^{144t-6} V(1)
\xrightarrow{v_2^{9t}} \Sigma^{-6}V(1) \xrightarrow{\nu} S^0
\\
\beta_{9t+1} & : S^{144t+10} \xrightarrow{v_2} \Sigma^{144t-6} V(1)
\xrightarrow{v_2^{9t}} \Sigma^{-6}V(1) \xrightarrow{\nu} S^0
\\
\beta_{9t+2} & : S^{144t+26} \xrightarrow{v_2} \Sigma^{144t+10} V(1)
\xrightarrow{v_2^{9t}} \Sigma^{10}V(1) \xrightarrow{v_2^*} S^0
\\
\beta_{9t+5} & : S^{144t+74} \xrightarrow{v_2^5} \Sigma^{144t-6} V(1)
\xrightarrow{v_2^{9t}} \Sigma^{-6}V(1) \xrightarrow{\nu} S^0
\\
\beta_{9t+6} & : S^{144t+90} \xrightarrow{v_2^5} \Sigma^{144t+10} V(1)
\xrightarrow{v_2^{9t}} \Sigma^{10}V(1) \xrightarrow{v_2^*} S^0
\end{align*}
\end{proof}

It should be the case that the elements $\beta_{9t+3}$ exist, but we are 
unable to deduce this from the existence of our self
map on $V(1)$.  Oka \cite{Oka} indicates that if the complex $M(3,v_1^3)$ has a
$v_2^9$ self map then the elements $\beta_{9t+3}$ exist.
Shimomura's computations of $\pi_*(L_2V(0))$ \cite{Shimomura:V(0)} 
demonstrate that the elements
$\beta_{9t+3}$ are present in $\pi_*(L_2(S^0))$.
Shimomura \cite{Shimomura:V(1)} proves that $\beta_i$ cannot be a
permanent cycle for $i \equiv 4,7,8 \pmod 9$.  In \cite{Shimomura:beta_1},
many relations amongst the $\beta_i$'s are investigated contingent on the
existence of the self-map constructed in this paper.  In particular,
Shimomura proves that if the elements $\beta_i \beta_1^4$ are permanent
cycles for $i \equiv 1 \pmod 9$, then they are non-trivial.  These elements
are permanent cycles by Corollary~$\ref{cor:betas}$.  
They should be regarded as the
substitutes for the the Adams-Novikov elements $\beta_{9t+4}$, 
which fail to exist.

Some remarks as to how this paper came to be written are in order.
The main result of this paper was the subject of the second author's
dissertation completed at Northwestern University under the direction of
Mark Mahowald.  The first author required the result for his dissertation
work at the University of Chicago under the direction of J. Peter May. 
Certain errors and
gaps in the original work needed to be corrected.  In the original thesis, 
the second author's main
technique for obtaining differentials in the ASS was to lift differentials
from the ANSS using a technical lemma called the `ladder lemma'.  We were
unable to make the proof of this lemma rigorous, and so the product rule
(Theorem~\ref{ProductRule}) is used instead for many of the
differentials.

The authors would like to thank Mark Mahowald for his constant
encouragement and assistance in this project.  We also thank Katsumi
Shimomura, for his useful correspondence concerning the construction of the
$\beta_i$'s from the self map.  The first author would also like to thank his
advisor, J. Peter May, for many useful conversations related to this
project.  We are also appreciative of many useful comments made by the
referee, including pointing out a substantially simpler proof of
Lemma~\ref{lem:dv_2^3b_1}.

\begin{conventions}
Throughout this paper we shall always be working in the stable homotopy
category localized at the prime $3$,
and all homology will be with $\FF_3$ coefficients.  We shall use the
following abbreviations.
\begin{verse}
ASS: Adams spectral sequence \\
ANSS: Adams-Novikov spectral sequence \\
BSS: Bockstein spectral sequence \\
AHSS: Atiyah-Hirzebruch spectral sequences
\end{verse}
The dual Steenrod algebra will be denoted by $A_*$.
If $X$ is a spectrum, we will often use the notation $\ext(X)$ to
represent $\ext_{A_*}(\FF_3, H_*(X))$, the $E_2$ term of the ASS for
computing $\pi_*(X)$.  We will denote the $E_r$ term of this ASS by
$E_r(X)$.
We shall use the notation $\doteq$ to indicate two quantities are equal up
to multiplication by a unit in $\FF_3$.  

Finally, in Section~\ref{sec:E_2} we give many elements in $H^*(P_*)$ names
which are derived from their May spectral sequence names.  In all but one
case, the sign of the element corresponding to a name coincides with the
element whose Curtis algorithm representative has a leading term with a
positive sign.  The one exception we make is the element called $k_0$ in
bidegree $(2,20)$.  We work under the sign convention that $k_0$ is
detected by the lambda algebra element $-\lambda_4 \lambda_3$.  The reason
we make this exception is so that certain relations
(Equation~\ref{eq:relations}) are more uniform.
\end{conventions}

\tableofcontents

\section{The Adams spectral sequence of $eo_2 \wedge V(1)$}\label{sec:eo_2}

In this section we will define $eo_2 \wedge V(1)$ and compute its ASS.  
The method of computing the $E_2$ term of the ASS is to produce
a finite complex $Y(2)$ such that upon 
smashing it with $eo_2 \wedge V(1)$ we get a wedge
of $k(2)$'s.  
We know the ASS of this object, and we recover the
$E_2$ term of the ASS for $\pi_*(eo_2 \wedge V(1))$ 
by forming a periodic resolution of $S^0$
out of copies of $Y(2)$.	
There is nothing original in this section.  Most of the material here was
first discovered by Hopkins, Mahowald, and Miller, but remains unpublished.

\begin{figure}
\begin{center}
  \includegraphics[width=4.9in]{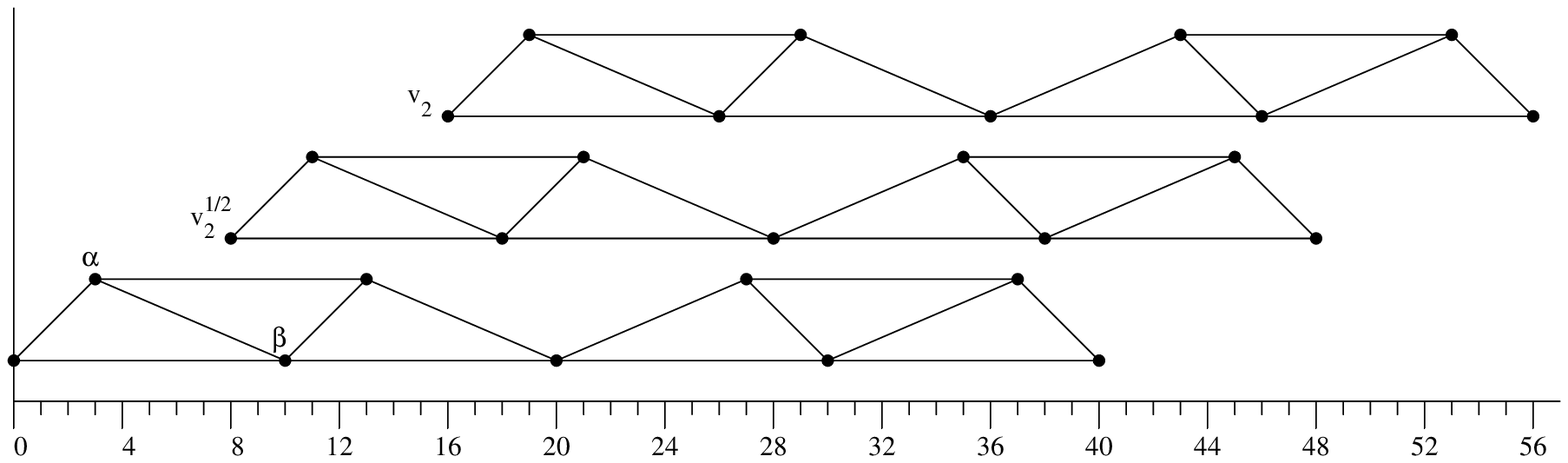}
\end{center}
\caption{${EO_2}_*(V(1))$ is generated as a free module over
$\FF_3[v_2^{\pm 9/2}]$ on the above pattern.} 
\label{fig:EO_2^V(1)}
\end{figure}

Let $E_2$ be the Hopkins-Miller spectrum at $p=3$.  It represents 
a Landweber exact
cohomology theory whose coefficient ring is 
$$ {E_2}_* = \WW_{\FF_9}[[u_1]][u,u^{-1}] $$
where $\abs{u_1} = 0$ and $\abs{u} = -2$.  
Here $\WW_{\FF_9}$ is the Witt ring with residue field $\FF_9$.
Fix a primitive $8^\mathrm{th}$ root of unity $\br{\omega}$ in $\FF_9$.  
We will refer to it's Teichm\"uller lift in $\WW_{\FF_9}$
also as $\omega$.
The element $\br{\omega}$ satisfies the relation
$$ \br{\omega}^2 + \br{\omega} + 2 = 0$$
in $\FF_9$.
The spectrum $E_2$ is a $BP$-ring spectrum, and
the map $\Phi:BP \rightarrow E_2$ has the following effect on coefficient
rings.
\begin{align*}
\Phi(v_1) & = u^{-2}u_1 \\
\Phi(v_2) & = u^{-8} \\
\Phi(v_i) & = 0, \quad \text{for $i > 2$}
\end{align*}
Let $\MS_2$ be the Morava stabilizer group.  It is the automorphism group
of the Honda height $2$ formal group law $F_2$ over $\FF_9$, and is
contained in the non-commutative algebra
$$ \WW_{\FF_9}\langle S \rangle /(S^2 = p, Sa = \sigma(a)S) $$
as the multiplicative group of units.  Here $\sigma$ is the a lift of the
Frobenius map.
The Galois group $Gal = Gal(\FF_9/\FF_3)$ acts on $\MS_2$ by acting on
$\WW_{\FF_9}$.  It is cyclic of order $2$ generated by the Frobenius
automorphism $\sigma$.  One may form the semi-direct product
$$ \GG_2 = \MS_2 \rtimes Gal. $$
The spectrum $E_2$ is an $E_\infty$ ring spectrum, and the group 
$\GG_2$ acts on $E_2$ via $E_\infty$
maps.  The spectrum $E_2$ and the $A_\infty$ action of 
$\MS_2$ are presented in
\cite{Rezk}.  There is a maximal finite 
subgroup $G_{12} < \MS_2$ of order $12$ which is
isomorphic to $C_3 \rtimes C_4$.  It is generated by an element $s$ of
order $3$, and an element $t$ of order $4$, given by the following
formulas.
\begin{align*}
s & = - \frac{1}{2}(1+\omega S) \\
t & = \omega^2
\end{align*}
These elements, as well as
specific formulas for their action on ${E_2}_*$, are given in
\cite{GoerssHennMahowald}.  
The subgroup $G_{12}$ is not invariant under the Galois action, so 
following \cite{GoerssHennMahowaldRezk}
we
instead investigate a maximal finite subgroup $G_{24} < \GG_2$ (of order
$24$) which
contains $G_{12}$ and fits into the following (non-split) 
short exact sequence.
$$ 1 \rightarrow G_{12} \rightarrow G_{24} \rightarrow Gal \rightarrow 1 $$
The subgroup $G_{24}$ is generated by the elements $s$, $t$, and $\psi$,
where we define
$$ \psi = \omega \sigma \in \GG_2. $$
The spectrum $EO_2$ is defined to be the homotopy fixed
point spectrum $E_2^{hG_{24}}$.  
A complete computation of the homotopy of $EO_2$, and its ANSS,
is given in \cite{GoerssHennMahowaldRezk}.   
In \cite{GoerssHennMahowald}, Goerss, Henn, and Mahowald compute the ANSS for
${EO_2}_*(V(1))$, but their approach must be
modified since they use $G_{12}$ instead of $G_{24}$.  
Their results may be conveniently summarized in
Figure~\ref{fig:EO_2^V(1)}.
In this figure, 
dots represent additive $\FF_3$ generators,
lines of length $3$ represent multiplication by $\alpha_1$,
lines of length $7$ represent the Toda Bracket $\langle \alpha_1, \alpha_1,
- \rangle$, and lines of length $10$ represent multiplication by $\beta_1$.
These products are given by the $S$-module structure.  The homotopy is
periodic with periodicity generator $v_2^{\pm 9/2}$ of degree $72$ on the
displayed pattern.

We want an Adams spectral sequence, but unfortunately the $v_2$-periodicity
in $EO_2$ makes its homology trivial.  We therefore need to take a
connective cover.  There is a nice connective cover of $EO_2$ called
$eo_2$ which has been constructed by Hopkins, Mahowald, and others.
Since the details of this construction are quite involved, 
we will instead
\emph{define} $eo_2 \wedge V(1)$ to be the connective cover of $EO_2
\wedge V(1)$.  

Since there is a gap in the homotopy of $EO_2\wedge
V(1)$ between the $56$ stem and the $72$ stem (and hence by the periodicity
of $EO_2$'s homotopy groups there is a gap between the $-16$ stem and $0$
stem), taking the connective cover
removes the periodic copies of the homotopy in negative dimensions.  
We
remark that the
reason that we cannot just define $eo_2$ to be the connective cover of
$EO_2$ is that the there are infinitely many copies of $BP \langle 1
\rangle $ in the
homotopy supported on negative periodicity generators whose unwanted homotopy
eventually appears in positive degrees.  Smashing with $V(1)$
kills all of this troublesome $v_1$-periodic homotopy.

We will now produce a finite complex $Y(2)$ which, when smashed with
$EO_2$, splits as a wedge of Morava $K$-theories.  
Let $\gamma: S^4 \rightarrow BO$ be a
generator of $\pi_4(BO) = \ZZ$.  The map $\gamma$ extends to a loop map
$\gamma_\infty:
\Omega S^5 \rightarrow BO$.  Let $J_i(S^4) \hookrightarrow \Omega S^5$ be
the $i^{\mathrm{th}}$ filtration of the James construction \cite{James}.  
Then $\gamma$ restricts to a map $\gamma_i: J_i(S^4) \rightarrow BO$.
Let $Y(i)$ be the Thom spectrum $(J_i(S^4))^{\gamma_i}$.
Then the homology of the ring spectrum $Y(\infty)$ is given by
$$ H_*(Y(\infty)) = \FF_3[b_2]$$
where $b_2$ has degree $4$.
The homology $H_*(Y(i))$ is the additive subgroup generated by 
$b_2^k$ for $0 \le k \le
i$.

There are maps $Y(i) \wedge Y(j) \rightarrow Y(i+j)$ induced from the maps
$J_i(S^4) \times J_j(S^4) \rightarrow J_{i+j}(S^4)$.  The complex 
$Y(1)$ is just the
Thom spectrum $(S^4)^\gamma$, and it is a two cell complex whose top cell
is attached to the bottom cell by the attaching map which is the image of
$\gamma$ under the $J$-homomorphism.  Therefore, $Y(1) = S^0
\cup_{\alpha_1} e^4$.  It follows that the dual Steenrod operation $P^1_*$
acts on $H_*(Y(1))$ by the formula
$$ P^1_*(b_2) = 1. $$
Using the map $Y(1) \wedge Y(1) \rightarrow Y(2)$, we obtain the following
formulas for the dual action of the Steenrod algebra on $H_*(Y(2))$.
\begin{align*}
P^1_* (b_2) & = 1 \\
P^1_* (b_2^2) & = 2 b_2 \\
P^2_* (b_2^2) & = 1
\end{align*}
In particular, we have the CW decomposition 
$Y(2) = S^0 \cup_{\alpha_1} e^4 \cup_{2\alpha_1} e^8$.

Our interest in $Y(2)$ arises from the following proposition.

\begin{prop}
There is a splitting
$$ EO_2\wedge V(1) \wedge Y(2) \simeq K(2) \vee \Sigma^8 
K(2). $$
\end{prop}

\begin{proof}
One can easily compute $\pi_*(EO_2 \wedge V(1) \wedge Y(\infty))$ from
the AHSS arising from the cellular
filtration of $Y(\infty)$, but the associated graded arising from this 
filtration gives too much
ambiguity for our purposes.  We therefore will use instead the 
homotopy
fixed point spectral sequence
$$ H^*(G_{24}; \pi_*(E_2 \wedge V(1) \wedge Y(\infty))) \Rightarrow \pi_*(EO_2
\wedge V(1) \wedge Y(\infty)) $$
where 
$$\pi_*(E_2 \wedge V(1) \wedge Y(\infty)) = \FF_9[u, u^{-1}, b_2]. $$
In \cite{GoerssHennMahowaldRezk}, the action of $G_{24}$
on ${E_2}_*(V(1))$ is given by the following formulas.
\begin{align*}
s_*(u) & = u \\
t_*(u) & = \br{\omega}^2 u \\
\psi_*(u) & = \br{\omega} u
\end{align*}
The elements $s,t \in G_{12}$ correspond to the automorphisms
\begin{align*}
s(x) &= x +_{F_2} u^{-2} \br{\omega} x^3 = x + u^{-2} \br{\omega} x^3 + 
\mathcal{O}(x^4) \\
t(x) &= \br{\omega}^2 x
\end{align*}
of the Honda height $2$ formal group $F_2$ over
$\FF_9[u,u^{-1}]$, with $3$-series $[3]_{F_2} = u^{-8} x^9$.  
Under the canonical map of Thom spectra $Y(\infty)
\rightarrow MU$, $b_2$ maps to the element of the same name in
$$ (E_2 \wedge V(1))_*(MU) = \FF_9[u, u^{-1}][b_1, b_2, b_3, \ldots ] $$
where the generators $b_i$ coincide with those of Adams in \cite[II.4.5]{Adams}.
These $b_i$ correspond to the coefficient of $x^{i+1}$ in a \emph{strict} 
map of 
formal
groups, and as such, we have
\begin{align*}
s_*(b_2) & = b_2 + u^{-2} \br{\omega} \\
t_*(b_2) & = b_2 \\
\psi_*(b_2) & = b_2.
\end{align*}
Therefore, the fixed points are given by
$$ \pi_*(E_2 \wedge V(1) \wedge Y(\infty))^{G_{24}} =
\FF_3[(\br{\omega}^2 u^{-4})^{\pm 1}, b^3_2-\br{\omega}^2 b_2u^{-4}] 
\subset \FF_9 [u^{\pm1}, b_2]. $$
Define $a_4 = \br{\omega}^2 u^{-4}$ and $a_6 = b_2^3 - \br{\omega}^2 
b_2 u^{-4}$.
In \cite[1.4]{GoerssHennMahowald}, 
the $E_2$ term of the homotopy fixed point
spectral sequence for $EO_2 \wedge V(1)$ is computed to be
$$ H^*(G_{24}; \pi_*(E_2 \wedge V(1))) = 
\FF_3[a_4^{\pm 1}, \beta]\otimes E[\alpha] $$
where $\beta = \langle \alpha, \alpha, \alpha \rangle$.  
The cellular filtration of $Y(\infty)$ gives an Atiyah-Hirzebruch type
spectral sequence that allows one to compute $H^*(G_{24}; \pi_*(E_2 \wedge
V(1) \wedge Y(\infty)))$ from this.  The $d_4$'s in this spectral sequence
are multiplication by $\alpha$, and the $d_8$'s are given by the
application of the Massey product $\langle \alpha, \alpha, - \rangle$.
Thus we conclude that
$$ H^s(G_{24}; \pi_*(E_2 \wedge V(1) \wedge Y(\infty))) =
\begin{cases}
\pi_*(E_2 \wedge V(1) \wedge Y(\infty))^{G_{24}}, & s = 0 \\
0, & s > 0
\end{cases} $$
and 
$$\pi_*(EO_2 \wedge V(1) \wedge Y(\infty)) = \FF_3[a_4^{\pm 1},a_6]. $$
The spectrum $EO_2 \wedge Y(\infty)$ is a ring spectrum whose homotopy is
concentrated in Adams-Novikov filtration $0$, and since the obstruction
for $V(1)$ to be a ring spectrum lies in positive Adams-Novikov filtration
(\ref{prop:splitting}),
$EO_2 \wedge V(1) \wedge Y(\infty)$ is also a ring spectrum.  Its
homotopy is concentrated in even degrees, so it is complex-orientable
\cite{Adams}.  The complex orientation 
$$ \theta: BP \rightarrow E_2 \wedge V(1) \wedge Y(\infty) $$
for which $\theta_*(v_i) = 0$ for $i \ne 2$ and $\theta_*(v_2) = u^{-8}$ lifts
to a complex orientation
$$ \td{\theta}: BP \rightarrow EO_2 \wedge V(1) \wedge Y(\infty).$$
Here the effect on homotopy is given by $\td{\theta}_*(v_i) = 0$ and
$\td{\theta}_*(v_2) = -a_4^2 $.
There are maps (for $\epsilon = 0,1$)
$$ S^{12j + 8\epsilon} \wedge BP \xrightarrow{a_4^\epsilon a_6^j \wedge \phi}
EO_2 \wedge V(1) \wedge Y(\infty) $$
that extend to maps
$$ \psi_{12j + 8\epsilon}: \Sigma^{12j+8\epsilon}K(2) 
\rightarrow EO_2
\wedge V(1) \wedge Y(\infty). $$
These maps give a splitting
$$ EO_2 \wedge V(1) \wedge Y(\infty) = \bigvee_{j \ge 0} 
\Sigma^{12j}\left( K(2)
\vee \Sigma^8 K(2)\right). $$
The composite
$$ EO_2 \wedge V(1) \wedge Y(2) \rightarrow EO_2 \wedge V(1)
\wedge Y(\infty) \rightarrow K(2) \vee \Sigma^{8} K(2) $$
(the second arrow is projection onto the first two wedge summands) is an
equivalence.
\end{proof}

\begin{cor}\label{cor:eo_2Splitting}
There is a splitting
$$ eo_2 \wedge V(1) \wedge Y(2) \simeq k(2) \vee 
\Sigma^8 k(2).$$
\end{cor}

\begin{proof}
The spectrum $k(2) \vee \Sigma^8 k(2)$ 
is the connective cover of $K(2) \vee
\Sigma^8 K(2)$.  
The Atiyah-Hirzebruch spectral sequence for $(eo_2 \wedge
V(1))_*(Y(2))$ is easily computed, and one finds
$$ \pi_*(eo_2 \wedge V(1) \wedge Y(2)) = \FF_3[a_4]. $$
Therefore, $eo_2 \wedge V(1) \wedge Y(2)$ is the connective cover of
$EO_2 \wedge V(1) \wedge Y(2)$.  The previous proposition and the
uniqueness of the connective cover combine to give this corollary.
\end{proof}

\begin{rmk}
Hopkins and Miller, in \cite{HopkinsMiller}, prove the following stronger
result.
$$ eo_2 \wedge Y(2) \simeq {BP \langle 2 \rangle}^{\wedge}_{3}  
\vee \Sigma^8 {BP \langle 2
\rangle}^\wedge_{3} $$
\end{rmk}

We will now construct a resolution of the sphere spectrum out of
suspensions of $Y(2)$.
There are cofiber sequences
\begin{gather*}
S^0 \rightarrow Y(2) \rightarrow \Sigma^4 Y(1) \\
Y(1) \rightarrow Y(2) \rightarrow S^8 
\end{gather*}
where the first maps are the evident inclusions.
Splicing these together gives the following 
$2$-periodic resolution of $S^0$.
$$
\xymatrix{
S^0 \ar[d] & 
\Sigma^3 Y(1) \ar[l] \ar[d] & 
S^{10} \ar[l] \ar[d] & 
\Sigma^{13} Y(1) \ar[l] \ar[d] &
S^{20} \ar[l] \ar[d] &
\cdots \ar[l]
\\
Y(2) &
\Sigma^3 Y(2) &
\Sigma^{10} Y(2) &
\Sigma^{13} Y(2) &
\Sigma^{20} Y(2) &
}
$$
The homology long exact sequences associated to this resolution 
break up into short exact sequences as a result of the following lemma.

\begin{lem}
The maps
\begin{gather*}
H_*(eo_2 \wedge V(1)) \rightarrow H_*(eo_2 \wedge
V(1) \wedge Y(2)) \\
H_*(eo_2 \wedge V(1) \wedge Y(1)) \rightarrow H_*(eo_2 \wedge V(1) \wedge Y(2))
\end{gather*}
are injective.
\end{lem}

\begin{proof}
The natural map of Thom spectra $Y(\infty) \rightarrow MU$ makes $MU$ a
$Y(\infty)$-ring spectrum, and therefore the Eilenberg-MacLane spectrum
$H\FF_3 = H$ is a $Y(\infty)$-ring spectrum.  Thus there is a retraction
$$
\xymatrix{
H \ar[r]^-{1 \wedge \eta} \ar@/_2pc/[rrr]_{Id} 
& H \wedge Y(2) \ar[r] & H \wedge Y(\infty)
\ar[r]^-{\mu} & H
}$$
and we may conclude that $H_*(eo_2 \wedge V(1)) \rightarrow H_*(eo_2 \wedge
V(1) \wedge Y(2))$ is an inclusion.  The complex $Y(2)$ is, up to suspension,
Spanier-Whitehead self-dual, and so we also have that the projection 
map $H_*(eo_2 \wedge
V(1) \wedge Y(2)) \rightarrow H_*(\Sigma^8 eo_2 \wedge V(1))$ is
surjective, hence, the previous map in the cofiber sequence 
$$ H_*(eo_2 \wedge V(1) \wedge Y(1)) 
\rightarrow H_*(eo_2 \wedge V(1) \wedge Y(2))$$ 
must
be injective.
\end{proof}

We may therefore apply $\ext_{A_*}(\FF_3, H_*(eo_2 \wedge V(1) \wedge -))$
to this resolution, and get long exact sequences, hence a spectral
sequence.  Our spectral sequence takes the form (for $\epsilon = 0,1$)
\begin{equation}\label{eq:Y(2)ss}
 E_1^{2k+\epsilon,s,t} = \ext^{s,t}(eo_2 \wedge V(1) \wedge
\Sigma^{10k+3\epsilon} Y(2)) \Rightarrow 
\ext^{s+2k+\epsilon,t+2k+\epsilon}(eo_2 \wedge V(1)) 
\end{equation}
Applying Corollary~\ref{cor:eo_2Splitting}, and using the known
computation $\ext(k(2)) = \FF_3[v_2]$ 
where $\abs{v_2} = (1,17)$, we may
express the $E_1$ term of \ref{eq:Y(2)ss} by
$$ E_1^{*,*,*} = \FF_3[v_2, \beta] \otimes E[\alpha, a]. $$
The tridegrees of these elements are $\abs{v_2} = (0,1,17)$, $\abs{a} =
(0, 0, 8)$, $\abs{\beta} = (2,0,10)$, and $\abs{\alpha} = (1,0,3)$.  The
only possible differentials are
$$ d_1(\alpha \beta^i v_2^j a) \doteq \beta^{i+1} v_2^j $$
but this $d_1$ arises from the composite $Y(2) \rightarrow S^8 \rightarrow
\Sigma^8 Y(2)$.  The element $v_2^j a \in \ext(eo_2 \wedge V(1) \wedge Y(2))$ 
is born on the zero cell of $Y(2)$, and
so must map to zero when projected onto the $8$-cell of $Y(2)$.  Therefore,
the spectral sequence \ref{eq:Y(2)ss} collapses at $E_1$, and we are left
with a computation of the $E_2$-term of the ASS for $\pi_*(eo_2 \wedge
V(1))$.

\begin{figure}[tp]
\begin{center}
  \includegraphics[width=4.9in]{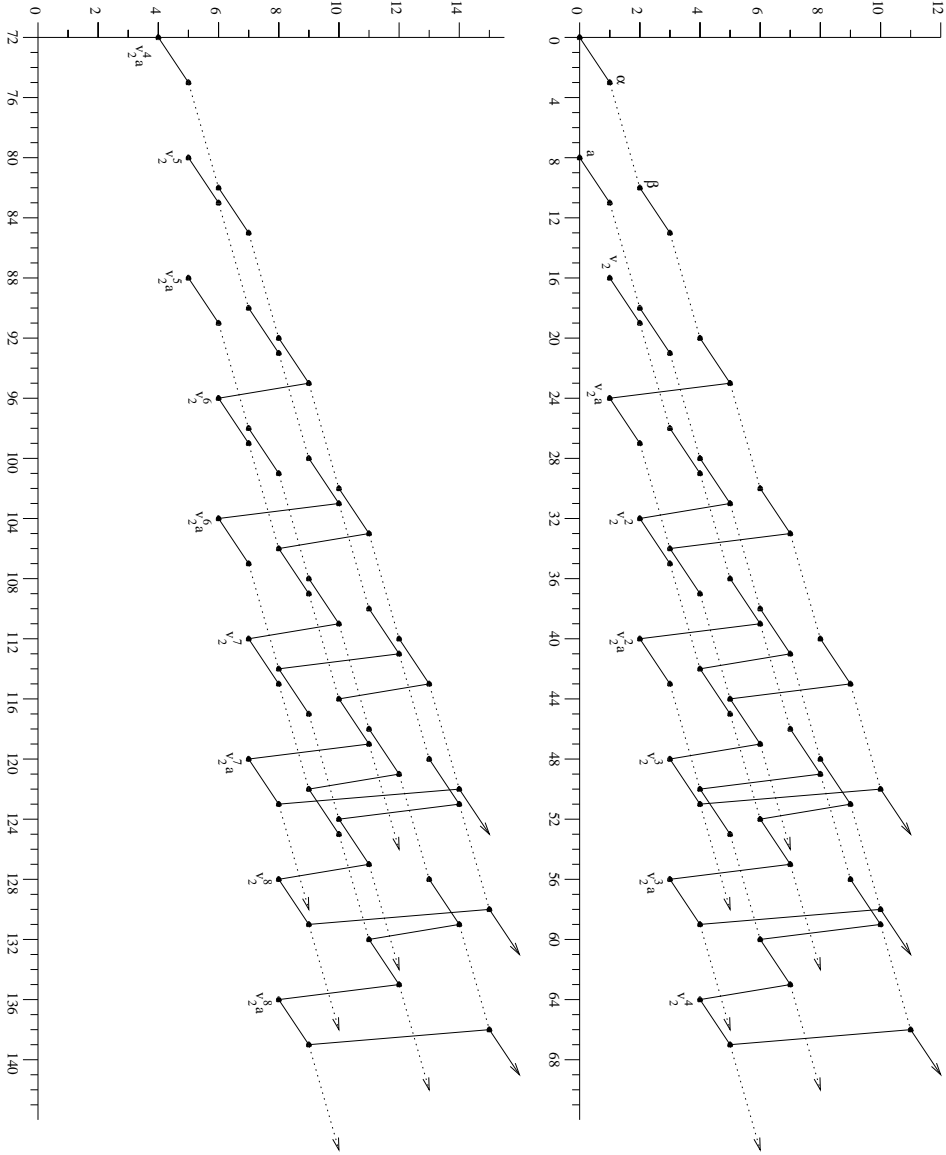}
\end{center}
\caption{The ASS for $\pi_*(eo_2 \wedge V(1))$ is a free module
over $\FF_3[v_2^9]$ on the displayed pattern.}
\label{fig:eo_2^V(1)ASS}
\end{figure}

The differentials in the ASS are
easily inferred from the differentials in the ANSS computed in
\cite{GoerssHennMahowald}.  Figure~\ref{fig:eo_2^V(1)ASS} displays the complete
ASS chart.
In this chart, $\FF_3$-generators are represented by dots, 
$\alpha$ multiplication is displayed with solid lines, and
the Massey product $\langle \alpha, \alpha, - \rangle$ is displayed with
dotted lines.  Solid lines of negative slope represent Adams differentials.
The $x$-axis represents the $t-s$ degree, and the $y$-axis represents the
homological degree $s$.

We remark that the ASS for $eo_2 \wedge V(1)$ is additively identical to
the ANSS.   The only difference is that $v_2$ has Adams filtration $1$,
whereas it has Adams-Novikov filtration $0$.

We finish this section with a computation of the effect of the $eo_2$ Hurewicz
homomorphism on the ASS.  We will use the
resolution of $S^0$ by the ring spectrum $Y(\infty)$.  
Mahowald, in \cite{Mahowald}, investigates a geometric Thom isomorphism
$$ Y(\infty) \wedge Y(\infty) \simeq Y(\infty) \wedge (\Omega S^5_+) $$
under which we may make the identification
$$ \pi_*(eo_2 \wedge V(1) \wedge Y(\infty) \wedge Y(\infty)) = \FF_3[a_4,
a_6, r] $$
Here $r$ has degree $4$.  We may regard $\FF_3[a_4,a_6, r]$ as being
contained in $\FF_3[u^{-1}, b_2]\otimes \FF_3[b_2]$ where $a_4$ and $a_6$ are
contained in the first factor as described earlier, and the element $r$
corresponds to $1 \otimes b_2$.  The main result of \cite{Mahowald} states
that the right unit of the associated Hopf algebroid
$$ \left( 
\FF_3[a_4,a_6], \FF_3[a_4,a_6,r] \right) $$
is given by the
following formulas where $b_2$ maps to $b_2 \otimes 1 + 1 \otimes b_2$.
$$
\xymatrix@R-2em{
\pi_*(eo_2 \wedge V(1) \wedge Y(\infty)) \ar[r]_-{\td{\Delta}} \ar@{=}[d]
& \pi_*(eo_2 \wedge V(1) \wedge Y(\infty) \wedge Y(\infty)) \ar@{=}[d] \\
\FF_3[a_4,a_6] & \FF_3[a_4,a_6,r] 
}$$
\begin{align*}
\td{\Delta}(a_4) = \Delta(\br{\omega}^2 u^{-4}) & = \br{\omega}^2 
u^{-4} \otimes 1 \\
 & = a_4 \\
\td{\Delta}(a_6) = \Delta(b_2^3 - \br{\omega}^2 b_2 u^{-4}) & = 
(b_2^3 - \br{\omega}^2 u^{-4} b_2) \otimes 1 
- \br{\omega}^2 u^{-4} \otimes b_2 + 1 \otimes b_2^3 \\
 & = a_6 - a_4r + r^3
\end{align*}
The Hurewicz image of $h_0$ is represented by $r$, and the Hurewicz image of 
$h_1$ is represented by 
$r^3$.
The $d_1$ supported by $a_6$ identifies the Hurewicz image of
$h_1$ with $h_0 a_4 = \alpha \cdot a$.  
This observation may be used to prove the
following proposition.

\begin{prop}\label{prop:eo_2Hurewicz}
The Hurewicz homomorphism
$$ h: \ext(V(1)) \rightarrow \ext(eo_2\wedge V(1)) $$
is described by 
\begin{alignat*}{3}
h(h_0) & = \alpha & \qquad h(b_0) & = \beta & \qquad h(v_2) & = v_2 \\
h(h_1) & = \alpha_1 a & \qquad h(g_0) & = \beta a
\end{alignat*}
\end{prop}

\section{Calculation of the 
Adams $E_2$ term $\ext_{A_*}(\FF_3, H_*(V(1)))$}\label{sec:E_2}

\begin{figure}[tp]
\begin{center}
  \includegraphics[height=8in]{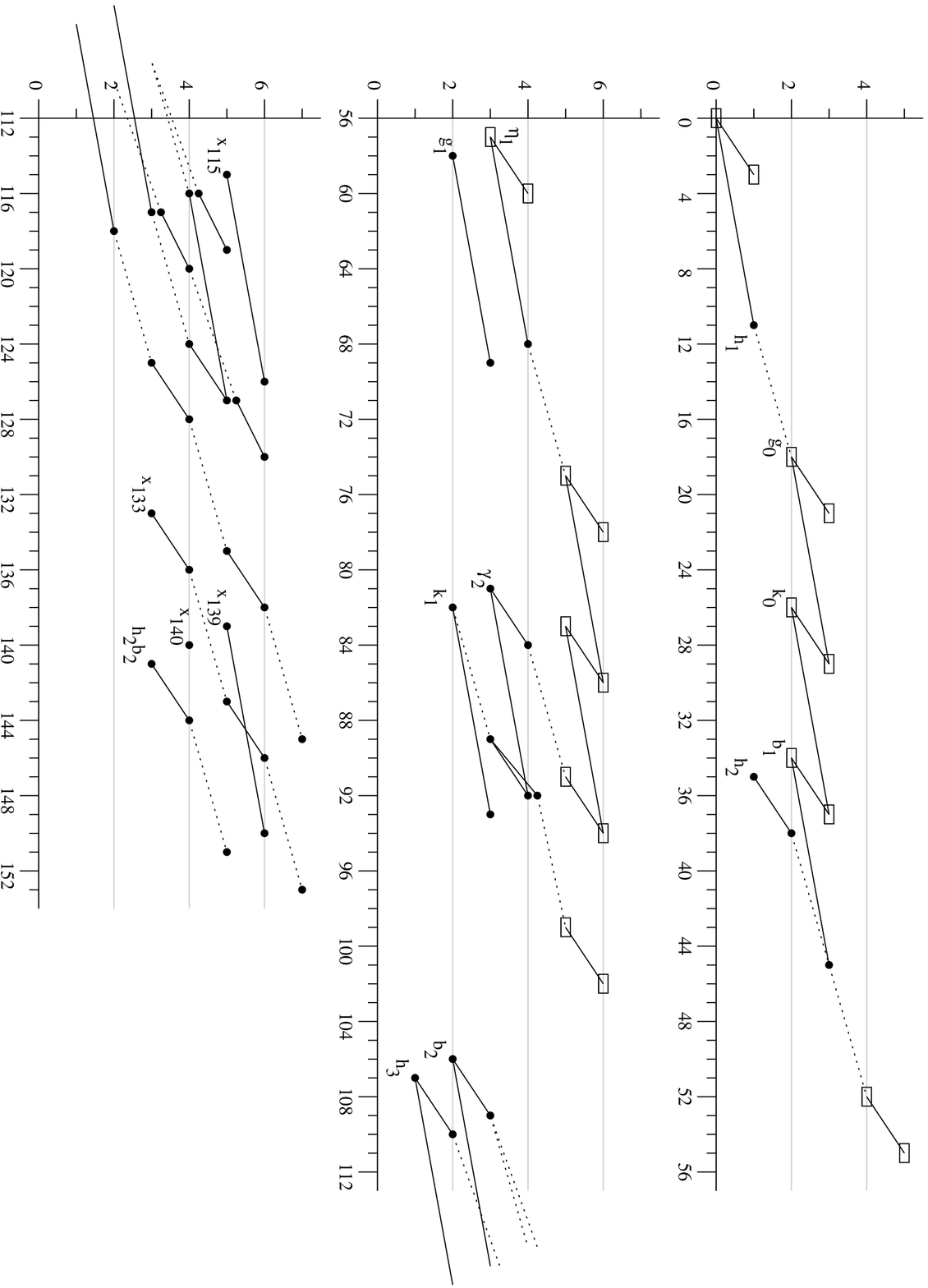}
\end{center}
\caption{$\ext_{P_*}(\FF_3, \FF_3)$, from Tangora's tables.}
\label{fig:extP}
\end{figure}

\begin{figure}[tp]
\begin{center}
  \includegraphics[height=8in]{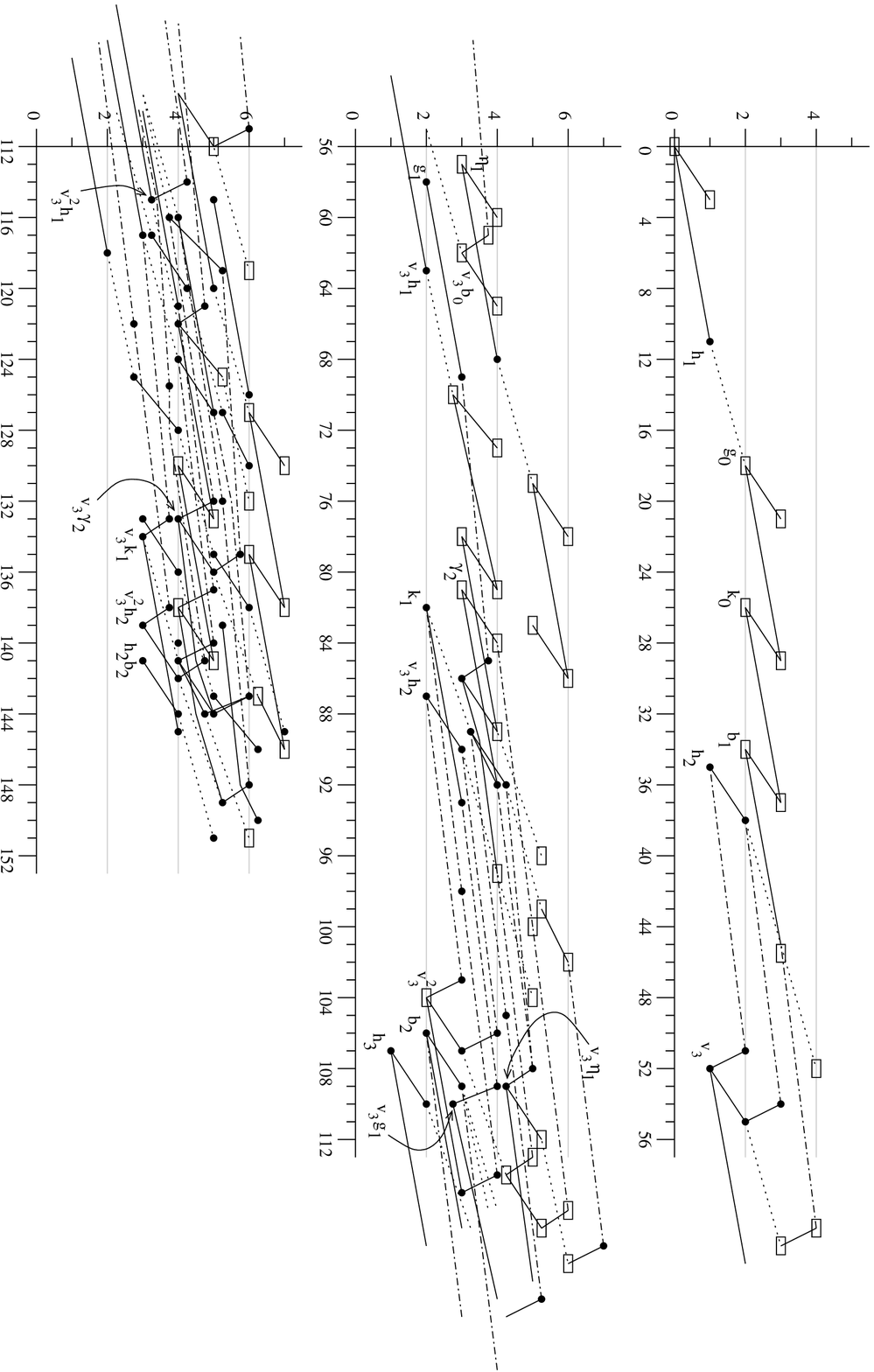}
\end{center}
\caption{The BSS $\ext_{A_*//E[\tau_0, \tau_1, \tau_2]}(\FF_3, \FF_3)
\otimes P[v_2] \Rightarrow \ext_{A_*//E[\tau_0, \tau_1]}(\FF_3, \FF_3)$.} 
\label{fig:bss}
\end{figure}

\begin{figure}[tp]
\begin{center}
  \includegraphics[height=8in]{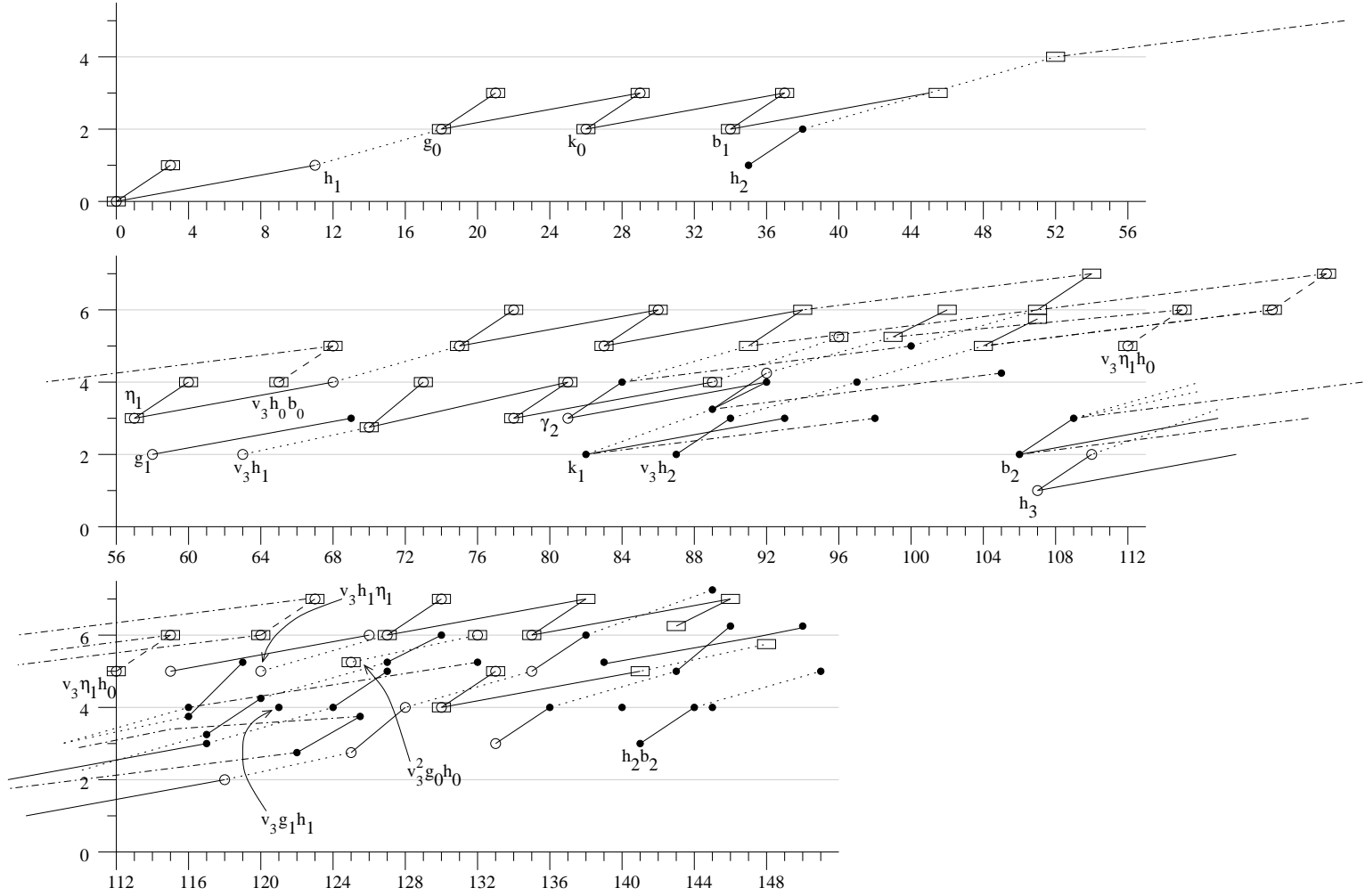}
\end{center}
\caption{The $E_2$ term $\ext_{A_*//E[\tau_0, \tau_1]}(\FF_3. \FF_3)$ of
the ASS for computing $\pi_*(V(1))$.}
\label{fig:assE2}
\end{figure}

The ASS for computing $\pi_*(V(1))$ has as its $E_2$ term
$$ \ext_{A_*}(\FF_3, H_*(V(1))) $$
As a
comodule over the dual Steenrod algebra, we have
$$ H_*(V(1)) = E[\tau_0, \tau_1]. $$
This is a subalgebra of the Steenrod algebra, but $V(1)$
is not a ring spectrum, so this algebra structure is not the consequence of
a geometric multiplication.  For the purposes of computing $\ext$, though,
we may use the algebra structure.  A change of rings isomorphism identifies
the $E_2$ term of the ASS as the cohomology of a Hopf algebra.
$$ \ext_{A_*}(\FF_3, E[\tau_0, \tau_1]) \cong \ext_{A_*//E[\tau_0,
\tau_1]}(\FF_3, \FF_3) = H^*(A_*//E[\tau_0, \tau_1]). $$
We may identify 
$$A_*//E[\tau_0, \tau_1] = P[\xi_1, \xi_2, \ldots]\otimes E[\tau_2, \tau_3,
\ldots].$$
The cohomology of this Hopf algebra is the cohomology of the subalgebra of the 
periodic
lambda algebra (see \cite{Gray} for a description of the periodic lambda
algebra) given by
$$ \br{\Lambda}_{(1)} = \langle \lambda_i, v_j \: : \: i \ge 0, j \ge 2
\rangle \subset \br{\Lambda}. $$
We only need the $E_2$ term of the Adams spectral sequence through the
dimension of $v_2^9$, which is $144$.  
Since the dimension of $v_4$ is $160$, 
the cohomology of $\br{\Lambda}_{(2)}$ coincides with $H^*(P_*) \otimes
P[v_3]$ in the range we are interested in.

Figure~\ref{fig:extP} displays
$H^*(P_*)$.  It was produced from Tangora's Curtis tables in
\cite{Tangora}.  The $y$ axis is the homological degree $s$, and the $x$-axis
is $t-s$, where $t$ is the internal degree.  
Solid lines of different slopes indicate 
multiplication by $h_0$ and $h_1$.
Dotted lines indicate the Massey product $\langle-, h_0, h_0\rangle$.
The element $b_0 = \langle h_0, h_0, h_0\rangle$, so $b_0$ multiplication
can be read off as composites of $h_0$ multiplication and application of
the above Massey product.  Rectangles indicate that a generator supports a
polynomial algebra on $b_0$, that is, all multiples of $b_0$ on the
generator are non-zero in the calculated range,
but they have not been written down in an
effort to make the chart less cluttered.   The Massey product 
representatives are the ones produced by using
the full tags produced by the Curtis algorithm.

In Figure~\ref{fig:extP}, certain generators have been given
names.  In our summary of conventions at the end of 
Section~\ref{sec:intro}, we
indicated that the signs of these elements will be chosen
(with the exception of $k_0$) so that the leading term of the corresponding
element of the lambda algebra has coefficient $+1$.  The following table
summarizes this choice of signs for some of the low dimensional generators,
by comparing our name, the lambda algebra name, and a Massey product
representation.
\vspace{.2in}

\begin{center}
\begin{tabular}{c|c|c}
\hline
Generator & Lambda Name & Massey Product \\
\hline \hline
$b_0$ & $\lambda_2\lambda_1+\lambda_1\lambda_2$ & $\langle h_0, h_0, h_0
\rangle$ \\
\hline
$g_0$ & $\lambda_2 \lambda_3$ & $\langle h_0, h_0, h_1 \rangle$ \\
\hline
$k_0$ & $-\lambda_4 \lambda_3$ & $\langle h_0, h_1, h_1 \rangle$ \\
\hline
$b_1$ & $\lambda_6 \lambda_3 + \lambda_3 \lambda_6$ & $\langle h_1, h_1,
h_1 \rangle $ \\
\hline
\end{tabular}
\end{center}
\vspace{.2in}

There are certain relations which may be read off of Figure~\ref{fig:extP}
up to sign.  We indicate the proper sign of some of the low dimensional
relations.  By choosing the sign of $k_0$ as we have, these relations look
more uniform.
\begin{equation}\label{eq:relations}
\begin{split}
h_1b_0 & = h_0g_0 \\
h_1g_0 & = h_0k_0 \\
h_1k_0 & = h_0b_1
\end{split}
\end{equation}

Figure~\ref{fig:bss} is a chart of the BSS
$$ \ext_{A_*//E[\tau_0, \tau_1, \tau_2]}(\FF_3, \FF_3)
\otimes P[v_2] \Rightarrow \ext_{A_*//E[\tau_0, \tau_1]}(\FF_3, \FF_3)
$$
It is straightforward to calculate the differentials of this spectral
sequence completely by explicitly finding the differentials in the periodic
lambda algebra $\br{\Lambda}_{(1)}$ and then finding the representatives in
the Curtis table.

In Figure~\ref{fig:bss}, the $E_1$-term consists of $H^*(A_*//E[\tau_0,
\tau_1, \tau_2])\otimes P[v_2]$ which is isomorphic to $H^*(P_*)\otimes
P[v_2, v_3]$ in our range of computation.  It is implicit in the chart that
every generator supports a $P[v_2]$, but all $v_2$ multiples are omitted
unless they are targets of differentials, or otherwise contribute to hidden
extensions.  When $v_2$ multiples are displayed, they are represented by
dash-dot lines.  Hidden extensions are represented by dashed lines, and
differentials are represented by negatively sloped solid lines.

We now have computed the $E_2$ term of the ASS.
It is displayed in Figure~\ref{fig:assE2}.  
Unlike in Figure~\ref{fig:bss}, in this chart $v_2$ is not implicit unless
specifically indicated.
Small solid dots on the chart
represent, like all of the previous charts, $\FF_3$-basis elements.  Small
circles represent polynomial algebras on $v_2$.  Otherwise $v_2$
multiplication is represented explicitly by dash-doted lines.  It should be
noted that if an element is represented by a circle, it does not mean that
that element supports infinitely many non-trivial multiplications by $v_2$. 
It just means that throughout the indicated range, all multiplications by
$v_2$ are non-trivial.

\section{The splitting of $D(V(1)) \wedge V(1)$}\label{sec:splitting}

The complex $V(1)$ may be visualized with the following cell diagram.
$$ \xymatrix@C-1em@R-2em{
0 & 1 &&& 5 & 6 \\
\circ \ar@{-}[r] & \circ \ar@{-}@/^1pc/[rrr] &&& \circ \ar@{-}[r] & \circ
}$$
Here the uncurved lines represent the Bockstein $\beta$ (attaching map
$\cdot 3$) and the curved line
represents the Steenrod operation $P^1$ (attaching map $\alpha_1$).
The top $V(0)$ is attached to the bottom $V(0)$ by $v_1$, but this is not
explicitly indicated in the cell diagram.
Let $D(V(1)) \simeq \Sigma^{-6}V(1)$ be the Spanier-Whitehead dual of
$V(1)$. 
In this section we will decompose $D(V(1)) \wedge V(1)$ into irreducible
subcomplexes.  Since $V(1)$ is self dual, we will have also provided a
splitting of $V(1) \wedge V(1)$.

\begin{figure}
\begin{center}
  \includegraphics[height=3in, clip, bb = 0 0 3.2in 3.6in]{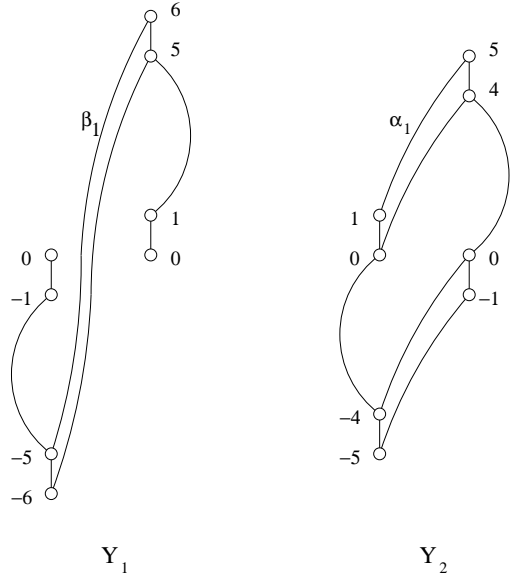}
\end{center}
\caption{The irreducible subcomplexes of $D(V(1)) \wedge V(1)$} 
\label{fig:splitting}
\end{figure}

Define finite complexes $Y_1$ and $Y_2$ as follows.
\begin{align*}
Y_1 & = \mathrm{cofiber}\left( \Sigma^{-1}V(1) \xrightarrow{\nu}
\Sigma^{4}V(0) \xrightarrow{\beta_1} \Sigma^{-6} V(1) \right) \\
Y_2 & = \mathrm{cofiber}\left( \Sigma^{-2}V(1) \xrightarrow{\alpha_1}
\Sigma^{-5} V(1) \right)
\end{align*}
Here $\nu$ is projection onto the top $V(0)$.
Figure~\ref{fig:splitting} displays cell diagrams of these complexes.

We will prove the following 
\begin{prop}\label{prop:splitting}
There is a splitting
$$ D(V(1)) \wedge V(1) \simeq Y_1 \vee Y_2. $$
\end{prop}

\begin{proof}
Since $D(V(1)) \wedge V(1) \simeq \Sigma^{-6} V(1) \wedge V(1)$, it
suffices to split the latter.
Consider the twist map
$$ \tau: V(1) \wedge V(1) \rightarrow V(1) \wedge V(1) $$
The map $\tau$ 
decomposes $V(1) \wedge V(1)$ into $+1$ and $-1$ eigenspaces $T_+$ and
$T_-$.  We claim that these are $\Sigma^{6}Y_1$ and $\Sigma^6 Y_2$, 
respectively.  More precisely, the self-maps
$(1+\tau)/2$ and $(1-\tau)/2$ are idempotents on $V(1) \wedge
V(1)$, and thus give a splitting 
$$ V(1) \wedge V(1) \simeq T_+ \vee T_-. $$
The structure of $T_-$ is straightforward from the action of the Steenrod
algebra.  To prove the existence of the $\beta_1$ attaching map in $T_+$,
we will use the secondary cohomology operation $\phi$ corresponding to the
Adem relation $P^2P^1 = 0$.  This secondary operation detects $\beta_1$.
In \cite{Thomas}, Thomas proves the Cartan formula
$$ \phi(xy) =
\phi(x)y+x\phi(y)+(P^1\beta(x))(\beta P^1 \beta(y))+(\beta P^1\beta(x))(P^1
\beta(y)). $$
Let $e_i$ denote the generator of $H^*(V(1))$ in dimension $i$.  Evaluating
$\phi$ on $e_0 \wedge e_0$, we get
$$ \phi(e_0\wedge e_0) = 0\wedge e_0+e_0\wedge 0+e_5\wedge e_6+e_6 \wedge
e_5 = e_5\wedge e_6 + e_6 \wedge e_5. $$
Similarly, we see $\phi(e_1 \wedge e_0 + e_0\wedge e_1) = -e_6 \wedge e_6$.
\end{proof}

\section{The product rule}\label{sec:ProductRule}

The statement of the product rule will require some notation related to
Adams resolutions which we will give presently.
Let $\br{H}$ be the fiber of the unit $\eta: S^0 \rightarrow H$, where $H$
is the Eilenberg-MacLane spectrum $H\FF_3$.
The standard
Adams resolution is defined by letting $W_s = \br{H}^{(s)}$ and then
defining $W_{s,r}$ to be the cofiber $\br{H}^{(s)}/\br{H}^{(s+r)}$.  
In particular, $W_{s,\infty} = W_s = \br{H}^{(s)}$ and $W_{s,1} = H \wedge
\br{H}^{(s)}$.
For a
spectrum $X$ the resolution may be written as
$$\xymatrix{
X \ar[d] & W_1 \wedge X \ar[d] \ar[l] & W_2 \wedge X
\ar[d] \ar[l] & \cdots \ar[l]
\\
W_{0,1}\wedge X & W_{1,1} \wedge X & W_{2,1} \wedge X
} $$
The ASS for $\pi_*(W_{s,r}\wedge X)$ is a truncated version of the ASS for
$\pi_*(X)$;  it is the spectral sequence obtained by only including the
$E_1^{i,j}$ terms for $s \le i < s+r$, and omitting differentials
supported on $E_r^{i,j}$ for $i < s$.
There are several important maps relating the spectra $W_{s,r}$.
\begin{gather*}
f_{s,r,k} : W_{s+r,k} \rightarrow W_{s,r+k} \\
g_{s,r,k} : W_{s,r+k} \rightarrow W_{s,r} \\
\partial_{s,r,k} : W_{s,r} \rightarrow \Sigma W_{s+r,k}\\
\mu_{s_1,s_2,r} : W_{s_1,r} \wedge W_{s_2,r} \rightarrow W_{s_1+s_2,r}
\end{gather*}
$\mu$ is induced by the product on $E$.  The remaining maps are compatible
in all of the ways one might expect them to be, and the sequence
$$ W_{s+r,k} \xrightarrow{f} W_{s,r+k} \xrightarrow{g} W_{s,r}
\xrightarrow{\partial} \Sigma W_{s+r,k} $$
is a cofiber sequence.

It is easy to see that an element $x \in \pi_*(W_{s,1} \wedge X)$ 
persists to the 
$E_r$ term of the ASS if and only if it lifts to an element $\tilde x \in
\pi_*(W_{s,r} \wedge X)$.  In fact, we have
$$
E_r^{s,t} = \frac{\im\{\pi_{t-s}(W_{s,r} \wedge X) 
\rightarrow \pi_{t-s}(W_{s,1} \wedge X)
\}}{\im\{\pi_{t-s+1}(W_{s-r,r} \wedge X)
\xrightarrow{\partial} \pi_{t-s}(W_{s} \wedge X) 
\rightarrow \pi_{t-s}(W_{s,1} \wedge X)
\}}
$$
and $d_r(x)$ is computed as $\partial (\tilde x) \in \pi_*(W_{s+r,1} \wedge
X)$.

Define, for $x_i$ in $\pi_*(V(1))$, 
$$ F(x_1,x_2) = \br{x_1} \cdot \brr{x_2} + (-1)^{\abs{x_1}} 
\brr{x_1} \cdot \br{x_2} 
\in \pi_*(V(0)) $$
where $\br{x_i}$ is the image of $x_i$ in $\pi_*(V(0))$ under
the projection $V(1) \rightarrow \Sigma^5 V(0)$, $\brr{x_i}$ is the image
of $\br{x_i}$ in $\pi_*(S)$ under the projection $V(0) \rightarrow
S^1$.  If $x_i \in \pi_*(W_{s_i,r} \wedge V(1))$, then $F(x_1,x_2)$ will be
regarded as an element of $\pi_*(W_{s_1+s_2,r} \wedge V(0))$.

\begin{thm}[{\bf Product Rule}]\label{ProductRule}
Suppose $x_i \in E_1(V(1))$ persist to the $E_r$-term of the ASS
and $\beta_1 \cdot F(x_1,x_2) = 0$, thought of an element of
$\pi_*(W_{s+2,r-2}
\wedge V(1))$ (the product is induced from the $V(0)$-module structure of
$V(1)$).  Then it follows that $\beta_1 \cdot F(x_1,x_2) \in
\pi_*(W_{s+2,r-1} \wedge V(1))$ lifts to an element $G(x_1,x_2)$ in 
$\pi_*(W_{s+r,1} \wedge V(1))$
and we have the following formula for $d_r(x_1 x_2)$.
$$ d_r(x_1 x_2) = (d_rx_1)\cdot x_2 + (-1)^{\abs{x_1}}x_1 \cdot (d_rx_2) -
G(x_1,x_2)
$$
\end{thm}

\begin{ex}\label{ex:dv_2^3}
We will use the product rule to compute $d_3(v_2^2)$.  The element
$v_2$ is a permanent
cycle for dimensional reasons.  
We have $\br{v_2} = \td{\beta_1} = h_1$ and $\brr{v_2} = \beta_1 =
b_0$.  Here we have given the ASS names of these elements.  Therefore,
$$ G(v_2,v_2) = (h_1 b_0 + b_0 h_1)b_0 = -h_1 b_0^2 $$
so the product rule says that
$$ d_3(v_2^2) = d_3(v_2)v_2 + v_2d_3(v_2) - G(v_2,v_2) = h_1 b_0^2. $$
One can also use the Hurewicz homomorphism 
$V(1) \rightarrow eo_2 \wedge V(1)$ to
get this formula.
\end{ex}

\begin{proof}[Proof of product rule]
Let $X$ be the cofiber of $\beta_1: \Sigma^{10}V(0) \rightarrow V(1)$.  
The closest substitute for a product on $V(1)$ is the map
$$ \mu: V(1) \wedge V(1) \rightarrow X $$
formed by projecting onto the wedge summand $\Sigma^6 Y_1$
(Proposition~\ref{prop:splitting}) and 
collapsing out the bottom two cells of the top $V(1)$.  
If $y_i$ are elements of $\pi_*(V(1))$, then the image
of $y_1 \wedge y_2 \in \pi_*(V(1) \wedge V(1))$ under the composition 
$$ V(1) \wedge V(1) \rightarrow X \rightarrow \Sigma^{11} V(0) $$
is $F(y_1,y_2)$.

We shall need various filtered forms of $X$.  Define $X_{s,r} =
W_{s,r}\wedge X$, and define
$\td{X}_{s,r}$ and $\tdd{X}_{s,r}$ to be the following cofibers.
\begin{gather*}
\Sigma^{10} W_{s,r} \wedge V(0) \xrightarrow{\beta_1} W_{s+2,r-2} \wedge
V(1) \rightarrow \td{X}_{s,r} \\
\Sigma^{10} W_{s,r} \wedge V(0) \xrightarrow{\beta_1} W_{s+2,r-1} \wedge
V(1) \rightarrow \tdd{X}_{s,r} \\
\end{gather*}
Note that the maps $\beta_1$ above may be chosen to raise the 
$s$-index by $2$ because they have Adams filtration $2$.
Then we have the following cofiber sequences (by Verdier's axiom).
\begin{gather*}
\td{X}_{s,r} \rightarrow X_{s,r} \rightarrow W_{s,2} \wedge V(1) \\
\tdd{X}_{s,r} \rightarrow \td{X}_{s,r} \rightarrow
  \Sigma W_{s+r,1}\wedge  V(1)
\\
\td{X}_{s,r+1} \rightarrow \td{X}_{s,r} \rightarrow \Sigma X_{s+r,1}
\\
\td{X}_{s,r+1} \rightarrow \tdd{X}_{s,r} \rightarrow \Sigma^{12} W_{s+r,1}
\wedge V(0) 
\end{gather*}

Since the Adams filtration of $\beta_1$ is greater than $0$, 
there is an equivalence $H
\wedge X \simeq H \wedge V(1) \vee H \wedge \Sigma^{11} V(0)$, thus
$X_{s,1}$ splits in a similar manner.  We need a splitting map
that behaves well with respect to the other maps floating about.
Consider the splitting $j$ induced on the cofibers below (the rows are
cofiber sequences).
$$
\xymatrix{
\td{X}_{s,r+1} \ar[r] \ar@{=}[d] &
\tdd{X}_{s,r} \ar[r] \ar[d] &
\Sigma^{12} W_{s+r,1} \wedge V(0) \ar@{.>}[d]_j \ar@/^3pc/[dd]^{Id}
\\
\td{X}_{s,r+1} \ar[r] &
\td{X}_{s,r} \ar[r] \ar[d] &
\Sigma X_{s+r,1} \ar[d]
\\
&
\Sigma^{12} W_{s+r,1} \wedge V(0) \ar@{=}[r] &
\Sigma^{12} W_{s+r,1} \wedge V(0)
} $$
We consider $j$ to be a nice splitting, because the projection $\nu$ it
induces onto the other wedge summand fills in the following 
triad of cofiber
sequences.
\begin{equation}\label{diag:Xtriad}
\xymatrix{
\td{X}_{s,r+1} \ar[r] \ar@{=}[d] &
\tdd{X}_{s,r} \ar[r] \ar[d] &
\Sigma^{12} W_{s+r,1} \wedge V(0) \ar@{.>}[d]_j 
\\
\td{X}_{s,r+1} \ar[r] \ar[d]&
\td{X}_{s,r} \ar[r]_{\partial_1} \ar[d]_{\partial_2} &
\Sigma X_{s+r,1} \ar@{.>}[d]_{\nu}
\\
\ast \ar[r] &
\Sigma W_{s+r,1} \wedge V(1) \ar@{=}[r] &
\Sigma W_{s+r,1} \wedge V(1)
}
\end{equation}

Let $\td{x}_i \in \pi_*(W_{s_i,r})$ be lifts of $x_i \in \pi_*(W_{s_i,1})$.
Let $s = s_1+s_2$.  It is useful to keep in mind the following diagram.
\vspace{.1in}

{\footnotesize
$$
\xymatrix{
W_{s,1}\wedge V(1) \ar[d]^\iota  &
W_{s,r}\wedge V(1) \ar[l]^g \ar[r]_\partial \ar[d]_\iota &
\Sigma W_{s+r,1}\wedge V(1) \ar[d]_\iota
\\
X_{s,1} \ar@/^1pc/[u]^\nu &
X_{s,r} \ar[l]^g \ar[r]_\partial &
\Sigma X_{s+r,1} \ar@/_1pc/[u]_\nu
\\
W_{s_1,1} \wedge W_{s_2,1} \wedge V(1)^{(2)} \ar[u]_\mu &
W_{s_1,r} \wedge W_{s_2,r} \wedge V(1)^{(2)} \ar[u]_\mu \ar[l]^g
\ar[r]_-{\partial} &
{\begin{array}{r}
 \Sigma W_{s_1+r,1} \wedge W_{s_2,r} \wedge V(1)^{(2)} \\ 
\vee \Sigma W_{s_1,r} \wedge W_{s_2+r,1} \wedge V(1)^{(2)}\ar[u]_\mu 
\end{array}}
} $$}
\vspace{.1in}

The element 
$x_1 \wedge x_2 \in \pi_*(W_{s_1,1} \wedge W_{s_2,1} \wedge V(1)^{(2)})$
lifts to $\td{x_1} \wedge \td{x_2} \in \pi_*(
W_{s_1,r} \wedge W_{s_2,r} \wedge V(1)^{(2)})$.  We then have
$$ \partial(\td{x_1} \wedge \td{x_2}) = d_r(x_1)\wedge x_2 +
(-1)^{\abs{x_1}}x_1 \wedge d_r(x_2). $$
The element $x_1 \cdot x_2$ is equal to $\nu\circ \mu (x_1 \wedge x_2) \in 
\pi_*(W_{s,1}\wedge V(1))$.  We want to
compute $d_r(x_1 \cdot x_2)$, which means we first need to lift $x_1 \cdot
x_2$ to $\pi_*(W_{s,r}\wedge V(1))$.  Now $\mu(\td{x_1} \wedge \td{x_2})$
is a lift of $\mu(x_1 \wedge x_2)$, but this element will not lift to a
lift of $x_1 \cdot x_2$ without a little modification.  
The following sequence is exact.
$$ \pi_*(\td{X}_{s,r}) \rightarrow \pi_{\ast -11}(W_{s,r} \wedge V(0))
\xrightarrow{\cdot \beta_1} \pi_{\ast - 1}(W_{s+2,r-2} \wedge V(1)) $$
Our assumption that $\beta_1 \cdot F(\td{x_1},\td{x_2}) \in
\pi_*(W_{s+2,r-2} \wedge V(1))$ is trivial implies that
$F(\td{x_1},\td{x_2})$ lifts to an element $\td{F} \in
\pi_*(\td{X}_{s,r})$.  Let $y$ be the image of $\td{F}$ in $X_{s,r}$, and 
define 
$$ z = \mu(\td{x_1} \wedge \td{x_2}) - y \in \pi_*(X_{s,r}). $$
We claim that (1) $z$ lifts to $\td{z} \in \pi_*(W_{s,r} \wedge V(1))$, 
and (2) $\td{z}$ is a
lift of $x_1 \cdot x_2 \in \pi_*(W_{s,1}\wedge V(1))$.

With regard to claim (1), we need only check that the image of $z$ in
$\pi_*(\Sigma^{11} W_{s,r} \wedge V(0))$ is zero.  The image of both
$\mu(\td{x_1} \wedge \td{x_2})$ and $y$ in $\pi_*(\Sigma^{11} W_{s,r} 
\wedge V(0))$ is $F(\td{x_1},\td{x_2})$, therefore the image of $z$, their
difference, is zero.  Claim (2) is established 
by noting that the sequence
$$ \pi_*(\td{X}_{s,r}) \rightarrow \pi_*(X_{s,r}) \rightarrow \pi_*(W_{s,2}
\wedge V(1)) $$
is exact.  Therefore the image of $y$ in $\pi_*(W_{s,2} \wedge V(1))$ 
is zero, so its image $\nu \circ g(y) \in \pi_*(W_{s,1} \wedge V(1))$ 
is zero.  
So, we have
$$ g(\td{z}) = \nu \circ \iota \circ g (\td{z}) = \nu \circ g(z) = \nu
\circ \mu(x_1 \wedge x_2) = x_1 \cdot x_2 $$
and claim (2) is established.

We are left with identifying $\partial (\td{z})$.  
We have
$$ \partial \td{z} = 
\nu \circ \partial (\mu(\td{x_1} \wedge \td{x_2}) - y)
= d_r(x_1)\cdot x_2 + (-1)^{\abs{x_1}}x_1 \cdot d_r(x_2) - \nu \circ
\partial (y) $$
We must evaluate $\nu \circ \partial(y)$.
In
Diagram~\ref{diag:Xtriad} the boundary maps $\partial_1$ and $\partial_2$
are displayed.  There is a map of cofiber sequences
relating $\partial_1$ to $\partial$ in the commutative diagram 
displayed below.
$$
\xymatrix{
\td{X}_{s,r} \ar[rr] \ar[dr]_{\partial_1} & & X_{s,r} \ar[dl]^\partial
\\
& \Sigma X_{s+r,1}
} $$
Therefore, $\partial y = \partial_1 \td{F}$.
Furthermore, Diagram~\ref{diag:Xtriad} reveals the relationship between
$\partial_1$ and $\partial_2$.  Thus we have $\nu \circ \partial_1 (\td{F}) =
\partial_2 (\td{F})$, and we just need an explicit description of the latter.
The map of cofiber sequences
$$
\xymatrix{
W_{s+r,1} \wedge V(1) \ar[r] \ar[d]_f & 
\tdd{X}_{s,r} \ar[r] \ar@{=}[d] &
\td{X}_{s,r} \ar[r]_{\partial_2} \ar[d] &
\Sigma W_{s+r,1} \wedge V(1) \ar[d]_f
\\
W_{s+2,r-1} \wedge V(1) \ar[r] &
\tdd{X}_{s,r} \ar[r] &
\Sigma^{11} W_{s,r} \wedge V(0) \ar[r]_{\beta_1} &
\Sigma W_{s+2, r-1} \wedge V(1)
} $$
tells us that $\partial_2 (\td{F})$ is a lift of $\beta_1 \cdot F(\td{x_1},
\td{x_2})$ to $\pi_*(W_{s+r,1}\wedge V(1))$ and as such, deserves to be called
$G(x_1,x_2)$.  This completes our verification of the formula.
\end{proof}

\begin{rmk}\label{rmk:ProductRule}
The theorem holds under a weaker assumption.  The proof of the theorem does
not require $x_1$ and $x_2$ to survive to $E_r$, but only that 
$\partial (x_1) \cdot x_2 + (-1)^{\abs{x_1}} x_1 \cdot \partial (x_2)$ have
Adams filtration greater than or equal to $s+r$.  We will need this
technical generalization for some of our applications of the product rule.
\end{rmk}

\section{Selected AHSS calculations of $\pi_*(V(1))$}\label{sec:AHSS}

In our calculation of differentials in the ASS it is helpful to know some
of the homotopy groups of $V(1)$.  The $3$-component of the homotopy groups
of spheres is known completely through the $108$ stem.  A table summarizing
these elements may be found in \cite[A3]{Ravenel}.  Thus one may write down
the $E_1$-term of the AHSS 
$$ E^1_{s,t} = \bigoplus_{\text{$s$-cells of V(1)}} \pi_{t+s}(S) \Rightarrow
\pi_{t+s}(V(1)) $$
in this range.  The complex $V(1)$ only has cells in dimensions $0$, $1$,
$5$, and $6$.  We shall denote an element in the $E^1$ term by the notation
$\gamma[k]$ where $\gamma \in \pi_*(S)$ and $k$ is the cell supporting it.
The differentials are determined by the attaching maps, and are
given by the following formulas.
\begin{align*}
d_1(\gamma[k]) & \doteq \begin{cases}
3\gamma[k-1] & k = 1,6 \\
0, & \text{otherwise}
\end{cases}
\\
d_4(\gamma[k]) & \doteq \begin{cases}
\alpha_1 \gamma [1] & k = 5 \\
0 & \text{otherwise}
\end{cases}
\\
d_5(\gamma[k]) & \doteq \begin{cases}
\langle \gamma, 3, \alpha_1 \rangle [1] & k = 6 \\
\langle \gamma, \alpha_1, 3 \rangle [0] & k = 5 \\
0 & \text{otherwise}
\end{cases}
\end{align*}

While a complete determination of the AHSS through the $108$ stem 
should be a relatively
straightforward task, we restrict ourselves to a few vicinities where we
need the data.  These partial charts are given on the next few pages, and 
are referred to in subsequent sections.

All but two of the differentials are immediate.  
We do not know if the dotted differential $(1)$ exists in
(\ref{AHSS:77-80})
because we are unsure
of whether or not $\beta_2^3 \in \pm\langle \beta_5, \alpha_1, 3 \rangle$.
We will see in the proof of Lemma~\ref{lem:stem139} that the 
differential $(1)$ must exist.
The only other differential
which isn't clear is $d_5(x_{68}[5])$ in (\ref{AHSS:68-73}).  We compute
$$
\alpha_1 \langle 3, \alpha_1, x_{68} \rangle  \doteq 3 \langle \alpha_1,
\alpha_1, x_{68} \rangle \ne 0. $$
This is a hidden extension in the ANSS for $\pi_*(S^0)$ in the computations
in \cite{Ravenel}.  The indeterminacy of $\langle 3, \alpha_1, x_{68}
\rangle$ is trivial, and the indeterminacy of $\langle \alpha_1, \alpha_1,
x_{68} \rangle$ is contained in $\alpha_1 \cdot \pi_{72}(S^0) = 3 
\cdot \pi_{75}(S^0)$, so
it doesn't enter into the above computation.  We conclude that $\langle 3,
\alpha_1, x_{68} \rangle \ne 0$, so it has no choice but to be a non-zero
multiple of $\beta_2^2 \beta_1^2$.
\vspace{.75in}

\begin{center}
{\bf Portions of the AHSS for $\pi_*(V(1))$}
\end{center}

\begin{equation}\label{AHSS:55-58}
\xymatrix@R-2em{
\underline{\text{Stem $55$}} & \underline{\text{Stem $56$}} & 
\underline{\text{Stem $57$}} 
& \underline{\text{Stem $58$}}
\\
\alpha_{14}[0] & \alpha_{14}[1] & \beta_2^2[5] \ar[dl] 
& \beta_2^2[6]
\\
\beta_2^2\alpha_1[0] & \beta_2^2\alpha_1[1] & \alpha_{13}[6] \ar[ul] \\
\beta_1^5[5] & \alpha_{13}[5] \ar[uul] & \\
\beta_2\beta_1^2\alpha_1[6] & \beta_1^5[6] 
}
\end{equation}

\begin{flushright}
{\it (cont'd on next page)}
\end{flushright}

\begin{landscape}
\begin{center}
{\bf Portions of the AHSS for $\pi_*(V(1))$, cont'd}
\end{center}
\vspace{.5in}

\begin{equation}\label{AHSS:63-68}
\xymatrix@R-2em{
\underline{\text{Stem $63$}} & \underline{\text{Stem $64$}} & 
\underline{\text{Stem $65$}} & \underline{\text{Stem $66$}}
& \underline{\text{Stem $67$}} & 
\underline{\text{Stem $68$}}
\\
\alpha_{16}[0] & \alpha_{16}[1] & \beta_2^2 \beta_1 \alpha_1[0] &
\beta_2^2\beta_1 \alpha_1[1] 
& \alpha_{17}[0] & x_{68}[0]
\\
\beta_2^2 \beta_1[1] & \alpha_{15/2}[5] \ar[ul] & 
\alpha_{15/2}[6] \ar[lu] \ar[l]^{\cdot 3} 
& & \beta_2^2 \beta_1[5] \ar[ul] & \alpha_{17}[1]
\\
&&&&& \alpha_{16}[5] \ar[uul]
\\
&&&&& \beta_2^2\beta_1[6]
}
\end{equation}
\vspace{.5in}

\begin{equation}\label{AHSS:68-73}
\xymatrix@R-2em{
\underline{\text{Stem $68$}} & \underline{\text{Stem $69$}} 
& \underline{\text{Stem $70$}} & \underline{\text{Stem $71$}} &
\underline{\text{Stem $72$}} & 
\underline{\text{Stem $73$}}
\\
x_{68}[0] & x_{68}[1] & \beta_2^2\beta_1\alpha_1[5] & \alpha_{18/3}[0] 
& \beta_2^2 \beta_1^2[0] & \beta_2^2 \beta_1^2[1]
\\
\alpha_{17}[1] & \alpha_{16}[6] \ar[l] & &
& \alpha_{18/3}[1] \ar[ul]_{\cdot 3} & x_{68}[5] \ar[ul]
\\
\alpha_{16}[5] & & & &
\alpha_{17}[5] \ar[uul] & \alpha_{17}[6] \ar[ul]
\\
\beta_2^2 \beta_1[6] &
}
\end{equation}

\begin{flushright}
{\it (cont'd on next page)}$\qquad \qquad \qquad$
\end{flushright}

\end{landscape}

\begin{center}
{\bf Portions of the AHSS for $\pi_*(V(1))$, cont'd}
\end{center}

\begin{equation}\label{AHSS:77-80}
\xymatrix@R-2em{
\underline{\text{Stem $77$}} & \underline{\text{Stem $78$}} 
& \underline{\text{Stem $79$}} & \underline{\text{Stem $80$}} \\
\beta_2^2\beta_1^2[5] & \beta_2^3[0] & \alpha_{20}[0] & \alpha_{20}[1] \\
\alpha_{18/3}[6] & \beta_2^2\beta_1^2[6] & \beta_2^3[1] & \alpha_{19}[5]
\ar[ul] \\
&& \beta_5[5] \ar@{.>}[uul]^{(1)} & \beta_2^2\beta_1^2\alpha_1[5] \\
&&& \frac{\beta_2^2\beta_1^2\alpha_1}{3}[5] \ar[uul] \\
&&& \beta_5[6]
}
\end{equation}
\vspace{.25in}

\begin{equation}\label{AHSS:87-90}
\xymatrix@R-2em{
\underline{\text{Stem $87$}} & \underline{\text{Stem $88$}} 
& \underline{\text{Stem $89$}} & \underline{\text{Stem $90$}} \\
\alpha_{22}[0] & \alpha_{22}[1] & \gamma_2\alpha_1[5] & \beta_6[0] \\
\beta_{6/2}[1] & \alpha_{21/2}[5] \ar[ul] & \beta_5\beta_1[5] &
\beta_{6/3}\alpha_1[5] \\
\beta_{6/3}[5] & \beta_{6/3}[6] \ar[ul] & \alpha_{21/2}[6] \ar[uul] \ar[ul]
& \langle \beta_2^3, \alpha_1, \alpha_1 \rangle [5] \\
\langle \beta_5, \alpha_1, \alpha_1 \rangle [6] &&& \gamma_2 \alpha_1[6] \\
\gamma_2[6] &&& \beta_5\beta_1[6]
}
\end{equation}
\vspace{.25in}

\begin{equation}\label{AHSS:97-100}
\xymatrix@R-2em{
\underline{\text{Stem $97$}} & \underline{\text{Stem $98$}} 
& \underline{\text{Stem $99$}} & \underline{\text{Stem $100$}} \\
x_{92}[5] & \beta_{6}\alpha_{1}[5] & \alpha_{25}[0] & \beta_2\beta_5[0] \\
\beta_{6/3}\beta_1[5] & \frac{\beta_6\alpha_1}{3}[5] & 
\langle x_{92},\alpha_1, \alpha_1 \rangle[0] & \alpha_{25}[1] \\
\alpha_{23}[6] & x_{92}[6] & \gamma_2\beta_1\alpha_1[5] & \langle x_{92}, 
\alpha_1, \alpha_1 \rangle[1] \\
\langle \beta_5\beta_1, \alpha_1, \alpha_1 \rangle[6] &
\beta_{6/3}\beta_1[6] & \beta_5\beta_1^2[5] & \alpha_{24/2}[5] 
\ar `l/16pt[u] `[uuul] [uuul] \\
\gamma_2\beta_1[6] && \beta_6\alpha_1[6] & \beta_{6/3}\beta_1\alpha_1[5] \\
&& \frac{\beta_6\alpha_1}{3}[6] \ar `l/16pt[u] `[uuuuul] [uuuuul]
& \gamma_2\beta_1\alpha_1[6] \\
&&& \beta_5\beta_1^2[6]
}
\end{equation}

\newpage

\section{The order of the $\beta_1$ action on $V(1)$}\label{sec:beta_1_order}

In this section we will prove the following proposition.

\begin{prop}\label{prop:beta_1_order}
The map 
$$ \beta_1^5:\Sigma^{50}V(1) \rightarrow V(1) $$
induced from smashing the map $\beta_1^5:S^{50} \rightarrow S^0$ with $V(1)$
is null.
\end{prop}

\begin{cor}\label{cor:beta_1_order}
Regarding $\pi_*(V(1))$ as a module over $\pi_*(S)$, we have the relation
$$ \beta_1^5 \cdot x = 0 $$
for all $x \in \pi_*(V(1))$.
\end{cor}

We remark that in $\pi_*(S)$ we have the relation $\beta_1^6 = 0$, and
$\beta_1^5$ is non-zero.  The power of $\beta_1^5$ in
Proposition~\ref{prop:beta_1_order} is minimal, since in $\pi_*(V(1))$ the
image of the element $\beta_1^4$ under the inclusion of the bottom cell is
non-trivial.

Corollary~\ref{cor:beta_1_order} follows from
Proposition~\ref{prop:beta_1_order} since the element $\beta_1^5\cdot x$
may be expressed as the following composite.
$$ S^{50+k} \xrightarrow{x} \Sigma^{50} V(1) \xrightarrow{\beta_1^5} V(1)
$$

We will first prove the following lemma.

\begin{lem}\label{lem:beta_1_order}
The element $\beta_1^5$ in $\pi_{50}(V(1))$ is trivial.
\end{lem}

\begin{proof}
There are no elements in the $50$ stem of Adams filtration greater than
$b_0^5$.  Therefore, it suffices to 
show that the element $b_0^5$ in the ASS for $\pi_*(V(1))$ is the
target of a differential.  In the ASS for $\pi_*(eo_2 \wedge V(1))$ there
is a differential
$$ d_6(v_2^3h_0) \doteq b_0^5. $$
Using the results of
Proposition~\ref{prop:eo_2Hurewicz}, we may conclude that if $v_2^3h_0$
supports no shorter differentials in the ASS for $V(1)$, then it must kill
$b_0^5$.  Upon investigating the $E_2$ term of the ASS for $V(1)$, we see
that there is no element in smaller Adams filtration that could be the target 
of a shorter differential.
\end{proof}

\begin{proof}[Proof of Proposition~\ref{prop:beta_1_order}]
We will demonstrate that the
Spanier-Whitehead adjoint of $\beta_1^5$
$$ \beta_1^5:S^{50} \rightarrow D(V(1)) \wedge V(1) $$
is null.  Let $X$ be the fiber of the composite
$$ V(1) \rightarrow \Sigma^5 V(0) \xrightarrow{\beta_1} \Sigma^{-5} V(0) $$
where the first arrow is projection onto the top $V(0)$.  By
Proposition~\ref{prop:splitting}, $X$ may be regarded as a subcomplex of
$Y_1$, which may in turn be regarded as a subcomplex of $D(V(1)) \wedge
V(1)$.
We wish to show that the composite
$$ S^{50} \xrightarrow{\beta_1^5} S^0 \hookrightarrow X \hookrightarrow Y_1
\hookrightarrow D(V(1)) \wedge V(1) $$ is null.  We will show that the
shorter composite $S^{50} \rightarrow Y_1$ is null.

Consider the following diagram, whose two bottom rows are cofiber
sequences.
$$\xymatrix@C+2em{
& S^{50} \ar@{.>}[ddl]_f \ar[d]_{\beta_1^5} \ar[ddr]^{*}
\\
& S^0 \ar[d] 
\\
\Sigma^{-6}V(0) \ar[r] \ar[d] & X \ar[r] \ar[d] & V(1)
\ar@{=}[d]
\\
\Sigma^{-6}V(1) \ar[r] & Y_1 \ar[r] & V(1)
}$$
In this diagram, the map $S^{50} \rightarrow V(1)$ is null by
Lemma~\ref{lem:beta_1_order}.  Therefore the lift $f$ exists making the
diagram commute.  We will complete the proof of the proposition once we
establish the following

\begin{claim*}
The image of the map $\pi_{56}(V(0)) \rightarrow \pi_{56}(V(1))$ is
trivial.
\end{claim*}

The claim follows easily from the AHSS for $\pi_*(V(0))$ and $\pi_*(V(1))$.  
A portion of the
AHSS for $\pi_*(V(0))$ is displayed below.
$$ \xymatrix@R-2em{
\underline{\text{Stem $55$}} & \underline{\text{Stem $56$}} & 
\underline{\text{Stem $57$}} \\
\alpha_{14}[0] & \alpha_{14}[1] & \\
\beta_2^2 \alpha_1[0] & \beta_2^2\alpha_1[1] 
}$$
There are no differentials, and $\pi_{56}(V(0))$ is of rank $2$.  We now
consider the image in $V(1)$. The same portion of the AHSS for $V(1)$ is
displayed in (\ref{AHSS:55-58}),
in which the same 
generators of $\pi_{56}(V(0))$ have been killed by differentials,
and the claim follows.
\end{proof}

\section{Proof that $v_2^9$ is a permanent cycle}\label{sec:mainthmpf1}

In this section we will prove that the element $v_2^9$ is a permanent cycle
in the ASS for $\pi_*(V(1))$.
We will let $h: V(1)
\rightarrow eo_2 \wedge V(1)$ be the Hurewicz homomorphism.  

We will first use the product
rule (\ref{ProductRule}) to determine $d_2(v_2^i)$.

\begin{lem}\label{lem:dv_2^i}
There are the following differentials on $v_2^i$ in the ASS.
\begin{alignat*}{3}
d_*(v_2) & = 0 
& \qquad d_3(v_2^2) & = h_1 b_0^2 
& \qquad d_2(v_2^3) & = -b_0 k_0 h_1 \\
d_2(v_2^4) & = -b_0 k_0 h_1 v_2
& \qquad d_2(v_2^5) & = -b_0 k_0 h_1 v_2^2
& \qquad d_2(v_2^6) & = b_0 k_0 h_1 v_2^3 \\
d_2 (v_2^7) & = b_0 k_0 h_1 v_2^4
& \qquad d_2(v_2^8) & = b_0 k_0 h_1 v_2^5 
& \qquad d_2(v_2^9) & = 0
\end{alignat*}
\end{lem}

\begin{proof}
These formulas are just obtained by iterated application of the product
rule.  The differential
$d_3(v_2^2)$ is computed in this manner in Example~\ref{ex:dv_2^3}. 
One then uses the following
formulas in the ANSS, which are derived in \cite[5.1.20]{Ravenel}
\begin{gather*}
\br{v_2^k} = \td{\beta_k} \equiv k h_1 v_2^{k-1} \pmod {v_1} \\
\brr{v_2^k} = \beta_k \equiv \binom{k}{2} v_2^{k-2} k_0 + kv_2^{k-1} b_0
\pmod {3,v_1}
\end{gather*}
to inductively determine $d_*(v_2^{k+1})$ from $d_*(v_2^k)$.  We should
point out that this formula differs by a sign from the formula in
\cite{Ravenel} because the elements we are referring to as $b_0$ and $k_0$
are normalized differently.
\end{proof}

\begin{rmk}\label{rmk:v_2^5}
In \cite{Oka}, Oka demonstrates that $v_2^5$ is a permanent cycle in the
ANSS for $\pi_*(V(1))$.  Above, we have shown that it supports a $d_2$ in
the ASS for $\pi_*(V(1))$.  In fact, there is a differential
$$ d_2(v_3 b_0 g_0) \doteq b_0 k_0 h_1 v_2^2 $$
and $v_2^5 \pm v_3 b_0 g_0$ is a permanent cycle in the ASS.  This
differential is established in Lemma~\ref{lem:v_2^5}.
\end{rmk}

We must eliminate the possibility that $v_2^9$ supports a $d_r$ for $r >2$.
We will make a list of all elements in the ASS in the $143$-stem of Adams
filtration greater than $11$.  It is given in the table below, with
references to the lemma that takes care of it, as well as the
Adams filtration (AF).
\vspace{.2in}

\begin{center}
\begin{tabular}{c|c|c}
\hline
AF & Element & Lemma \\
\hline \hline
29 & $h_0 b_0^{14}$ & \ref{lem:target1} \\
\hline
23 & $g_0 h_0 v_2^2 b_0^9$ & \ref{lem:target1} \\
\hline
22 & $b_1 h_0 v_2 b_0^9$ & \ref{lem:target2} \\
\hline
18 & $h_0 b_0^6 v_2^5$ & \ref{lem:target1} \\
   & $\eta_1 v_2 b_0^7$ & \ref{lem:target3} \\
   & $v_3 h_1 b_0^8$ & \ref{lem:target4} \\
\hline
17 & $k_0 h_0 v_2^4 b_0^5$ & \ref{lem:target5} \\
   & $\eta_1 k_0 b_0^6$ & \ref{lem:target6} \\
\hline
13 & $v_3 h_0 b_0^4 v_2^3$ & \ref{lem:target9} \\
\hline
12 & $g_0 h_0 b_0 v_2^7$ & \ref{lem:target1} \\
   & $\eta_1 g_0 v_2^3 b_0^2$ & \ref{lem:target8} \\
   & $v_3 k_0 h_0 v_2^2 b_0^3$ & \ref{lem:target7} \\
   \hline
\end{tabular}
\end{center}
\vspace{.2in}

\begin{lem}\label{lem:target1}
If $x \in E_2(V(1))$ is an element of the ASS for 
$V(1)$, and its Hurewicz image
$h(x) \in E_2(eo_2 \wedge V(1))$ is non-zero, then $x$ cannot be the
target of a differential supported by $v_2^9$.
\end{lem}

\begin{proof}
We have $h(d_r (v_2^9)) = d_r h(v_2^9)$, but $h(v_2^9) \in E_2(eo_2 \wedge
V(1))$ is a permanent cycle, so $d_r(v_2^9)$ must be in the Hurewicz kernel.
\end{proof}

\begin{lem}\label{lem:target2}
Suppose that $y \in E_6(V(1))$.  Then $h_0 b_0^3 y = 0$ in $E_6(V(1))$.
Similarly, if $z \in E_7(V(1))$, then $b_0^6 z = 0$ in $E_7(V(1))$.  In
particular, if $x$ is an element of 
$E_2(V(1))$ and $x = h_0 b_0^3y$ or $x= b_0^6 z$ for
$y$ or $z$ as before, than $x$ is not the target of a non-trivial $d_r$
for $r \ge 6$ or $r \ge 7$, respectively.
\end{lem}

\begin{proof}
The element $h_0 b_0^3$ dies in $E_6(S^0)$ (this is just the Toda
differential $d_5 b_1 \doteq h_0 b_0^3$).  Similarly, there is a
differential $d_6 (h_0 b_0 v_2^3) \doteq b_0^6$ giving the relation $b_0^6
= 0$ in $E_7(S^0)$.  Then use the $S^0$-module structure of $V(1)$.
\end{proof}

The following lemmas take care of the other possible targets.  We work from
highest to lowest Adams filtration to eliminate the possibility of
intervening differentials as we go along.

\begin{lem}\label{lem:target3}
In $E_5(V(1))$ there is a non-trivial differential
$$ d_5(\eta_1 v_2 b_0^7) \doteq k_0 v_2 b_0^{10}. $$
\end{lem}

\begin{proof}
In $E_5(S^0)$, we have $d_5 \eta_1 \doteq k_0 b_0^3$.  
The element $v_2 b_0^7$ is a
permanent cycle in the ASS for $V(1)$, so the differential follows from the
$S$-module structure of $V(1)$.
\end{proof}

\begin{lem}\label{lem:target4}
In $E_4(V(1))$ there is a non-trivial differential
$$ d_4(v_3 h_1 b_0^8) \doteq b_0^9 k_0^2. $$
\end{lem}

\begin{proof}
Our AHSS calculations (\ref{AHSS:68-73}) 
prove that $\pi_{72}(V(1)) = 0$.  
Therefore, the ASS for $\pi_*(V(1))$ should have no non-trivial permanent
cycles in the $72$-stem.  The $E_2$ term contains $k_0^2 b_0^2$, $b_0^3 k_0
v_2$, and $b_0^4 v_2^2$.  The element
$b_0^4 v_2^2$ supports a non-trivial $d_3$, and
$d_5 \eta_1 v_2 \doteq b_0^3 k_0 v_2$.  The only possibility for eliminating
$k_0^2 b_0^2$ is for $d_4(v_3 h_1 b_0) \doteq k_0^2 b_0^2$.  Therefore
$d_4(b_0^8 v_3 h_1) \doteq k_0^2 b_0^9$.
\end{proof}

For Lemma~\ref{lem:target5}
we need to know the differentials supported by $v_2^i g_0$ for small $i$.
These are given below.

\begin{lem}\label{lem:dg_0v_2^i}
We have the following Adams differentials on $g_0 v_2^i$ in $E_r(V(1))$.
\begin{alignat*}{3}
d_4(g_0 v_2) & = b_0^3 h_0 & \qquad 
d_3(g_0 v_2^2) & = b_0^2 k_0 h_0 & \qquad
d_2(g_0 v_2^3) & = -g_0 b_0 k_0 h_1 \\
d_3(g_0 v_2^4) & = 0 & \qquad
& & \qquad
d_3(g_0 v_2^6) & = v_2^4 b_0^2 k_0 h_0
\end{alignat*}
\end{lem}

\begin{rmk} 
The element $g_0v_2^4$
is actually a permanent cycle, 
and this should be regarded as anomalous.  
The AHSS element which it corresponds to is $\langle \beta_5, \alpha_1,
\alpha_1 \rangle[1]$, and this bracket is defined only because of 
the anomalous relation $\alpha_1 \beta_5 = 0$ in $\pi_*(S)$.
\end{rmk}

\begin{proof}
We will first explain how the term $G(g_0 v_2^i, v_2)$ and the 
term $G(g_0 v_2^i,
v_2^2)$ in the product rule is
computed.
The Adams-Novikov
element which detects $\br{g_0 v_2^i}$ is given by
$$ v_2^i b_0 h_0 - i v_2^{i-1} k_0 h_0 \pmod {v_1}. $$
(We will just work modulo $(v_1)$ since we will be mapping everything into
$V(1)$ for the product rule anyways.)
The Adams-Novikov element which detects $\brr{g_0 v_2^i}$ is given by 
$$ \binom{i}{2}v_2^{i-2} k_0 g_0 + i v_2^{i-1} b_0 g_0 \pmod {v_1}. $$
We recall from Lemma~\ref{lem:dv_2^i} the following formulas.
\begin{alignat*}{2}
\br{v_2} & = h_1 & \qquad
\br{v_2^2} & = -h_1 v_2 \\
\brr{v_2} & = b_0 & \qquad
\brr{v_2^2} & = k_0 - b_0 v_2
\end{alignat*}
Using the relations~\ref{eq:relations} and the relation $v_2b_1h_1=0$,
we may 
apply the product rule (\ref{ProductRule}) iteratively 
to get the requisite differentials.  Specifically, first
apply the product rule to $g_0\cdot v_2$, to get $d_4(v_2 g_0)$, then apply
the product rule to $(g_0v_2) \cdot v_2$ to get $d_3(g_0 v_2^2)$.  In
$E_2(S^0)$, $g_0$ supports a $d_2$, and $v_2^3$ supports a $d_2$ in
$E_2(V(1))$, thus $d_2(v_2^3 g_0)$ may be deduced from the $S$-module
pairing of Adams spectral sequences.  The problem is, we can no longer
apply the product rule to $v_2$ multiplication to get $d_3 v_2^4 g_0$.
However, we may instead apply the product rule to the product 
$v_2^2 \cdot v_2^2 g_0$ to get the
formula for $d_3(v_2^4 g_0)$ and similarly to $v_2^2 \cdot v_2^4 g_0$ to
get the formula for $d_3(v_2^6 g_0)$.
\end{proof}

In particular, we have the following lemma, which follows immediately.

\begin{lem}\label{lem:target5}
In $E_3(V(1))$ there is a differential
$$ d_3(v_2^6 b_0^3g_0) = k_0 h_0 v_2^4 b_0^5. $$
\end{lem}

\begin{lem}\label{lem:target6}
In $E_5(V(1))$ there is a differential
$$d_5(v_2^4 h_0 b_0^2 \eta_1) \doteq \eta_1 k_0 b_0^6.$$
\end{lem}

\begin{proof}
We need to compute the differential supported by $v_2^4 h_0$.  Observe that
the differential $d_6(v_2^3 h_0) \doteq b_0^5$ in the ASS for $eo_2 \wedge
V(1)$ lifts to a differential in
$E_6(V(1))$.  
The element $v_2^2 h_0$ must be a permanent cycle since there is nothing
for it to kill.  Hence, by the product rule, 
$$ b_0^5 \doteq d_6(v_2^3 h_0) = -G(v_2^2 h_0, v_2) = 
-b_0(\br{h_0 v_2^2} b_0 - \brr{h_0 v_2^2}h_1). $$
We conclude that the image of $\br{h_0 v_2^2}$ in $\pi_*(V(1))$ is
$\pm \beta_1^3$.  The image of $\brr{h_0 v_2^2}$ in $\pi_*(V(1))$ has to be in
Adams filtration greater than or equal to that of the image of $\br{h_0
v_2^2}$, and so we may conclude that the image of $\brr{h_0 v_2^2}$ is
actually zero.  Whereas $d_3(v_2^2)$ is non-zero, the differential
$(d_3(v_2^2)) v_2^2 h_0$ is
zero, and there are no permanent cycles in higher Adams filtration.
Therefore 
$(\partial(v_2^2)) \cdot v_2^2 h_0 = 0$, and we are in a position to
use the version of the product rule explained in
Remark~\ref{rmk:ProductRule}.  Since $\brr{v_2^2}$ is detected by $k_0$, we
have
$$ d_5(v_2^4 h_0) = -G(v_2^2, v_2^2 h_0) = -b_0(\br{v_2^2} \cdot \brr{v_2^2
h_0} + \brr{v_2^2} \cdot \br{v_2^2 h_0}) \doteq -b_0(k_0 \cdot b_0^3) = -k_0
b_0^4. $$
We then use the $S$-module structure of $V(1)$, and the differential $d_5
\eta_1 \doteq k_0 b_0^3$ to obtain
$$ d_5(\eta_1 \cdot v_2^4 h_0 b_0^2) = \pm k_0 h_0 b_0^5 v_2^4 \pm \eta_1 k_0
b_0^6. $$
By Lemma~\ref{lem:target5}, $k_0 h_0 b_0^5 v_2^4$ is the target
of a $d_3$.
\end{proof}

We will need a couple of lemmas to prove Lemmas~\ref{lem:target9} and
\ref{lem:target7}.

\begin{lem}\label{lem:dv_3h_0b_0}
In $E_4(V(1))$ there is a differential
$$ d_4(v_3 h_0 b_0) \doteq b_1 b_0^3. $$
\end{lem}

\begin{proof}
In the Adams spectral sequence for $\pi_*(V(1))$ in the vicinity of the 
$65$ stem, we have
the following elements and differentials.
$$\xymatrix@R-2em{
\underline{\text{Stem $63$}} & \underline{\text{Stem $64$}} & 
\underline{\text{Stem $65$}} & \underline{\text{Stem $66$}}
\\
b_0^6 h_0 & b_0^3 g_0 v_2 \ar[l] & b_0^3h_0v_2^2 & b_0^5 v_2 \\
b_0^2 v_2^2 h_1 & b_1 b_0^3 & b_0^2 h_0 k_0 v_2 & b_0^4 k_0 \\
b_1 b_0 h_0 v_2 & v_2^4 \ar[l] & b_0^2 b_1 h_1 & v_2^3 g_0 \ar[l] \\
v_3 h_1 & & v_3 h_0 b_0 \ar@{.>}[uul]^{(2)} & b_1 v_2^2
\ar@{.>}[uuul]_{(1)} \ar@{.>}[uul]
}$$
Here we know for dimensional reasons that $b_1$, $k_0$, $v_2$, $h_0 v_2^2$,
and $h_0k_0 v_2$ are
permanent cycles, and so all $b_0$ multiples of them are permanent cycles.
We have the solid differentials from Lemmas~\ref{lem:dv_2^i} and
\ref{lem:dg_0v_2^i}.  The only means by which the correct homotopy can be
achieved is for one of the dotted differentials $(1)$ to exist
and the dotted differential $(2)$ to exist.  The differential $(2)$ is the
desired differential.
\end{proof}

\begin{lem}\label{lem:dv_2^3b_1}
The element $v_2^3 b_1$ is a permanent cycle.
\end{lem}

\begin{proof}
In \cite{Ravenel} the element $\beta_{6/3}$ is shown to exist in 
$\pi_{82}(S^0)$.  It is detected modulo the ideal $(3,v_1)$ in the ANSS by
the element $\pm v_2^3b_1$.  Thus, the image of $\beta_{6/3}$ in 
$\pi_{82}(V(1))$
is detected in the ASS by the element $\pm v_2^3b_1$, and this element
must therefore be a permanent cycle.
\end{proof}

\begin{lem}\label{lem:target9}
In $E_2(V(1))$ there is a non-trivial differential
$$ d_2(v_3h_0b_0^4v_2^3) \doteq v_2h_0\eta_1k_0b_0^4. $$
\end{lem}

\begin{proof}
In Lemma~\ref{lem:dv_3h_0b_0}, we showed that $d_2(v_3h_0b_0) = 0$.  In
Lemma~\ref{lem:dv_2^i}, we showed that 
$$d_2(v_2^3) = -b_0k_0h_1 = -b_0b_1h_0. $$
We also know $\br{v_2^3}$, and hence $\brr{v_2^3}$, are trivial in $E_2$
(they are in higher Adams filtration).
We will now use the product rule (\ref{ProductRule}) to compute the
differential on $(v_3h_0b_0) \cdot v_2^3$.  
The term $G(v_3h_0b_0, v_2^3)$ is trivial in $E_2$ by the previous
considerations.
The following manipulations are
made possible with the hidden extension $h_0\cdot(v_3h_0b_0) \doteq
b_1g_0v_2$.
\begin{align*}
d_2((v_3h_0b_0)\cdot v_2^3) & = (v_3h_0b_0)\cdot(b_0b_1h_0) \\
& \doteq b_1^2g_0b_0v_2 \\
& \doteq v_2h_0\eta_1k_0b_0
\end{align*}
Now multiply by $b_0^3$.
\end{proof}

\begin{lem}\label{lem:target7}
In $E_5(V(1))$ there is a non-trivial differential
$$ d_5(v_3 k_0h_0v_2^2b_0^3) \doteq b_1 b_0^6 v_2^3. $$
\end{lem}

\begin{proof}
In $E_2(V(1))$, there is a Massey product
\begin{equation}\label{eq:MasseyProduct}
v_2^3b_1b_0^3 \in \langle b_1b_0^3, v_1, v_1^2 h_2 \rangle.
\end{equation}
We
are regarding the elements in the Massey product as detecting the
following maps in homotopy.
\begin{align*}
b_1b_0^3 & \leftrightarrow \beta_2\beta_1^3\alpha_1[5] 
: \Sigma^{64}V(0) \rightarrow
V(1) \\
v_1 & \leftrightarrow v_1:\Sigma^4V(0) \rightarrow V(0) \\
v_1^2h_2 & \leftrightarrow \td{\beta_3}: \Sigma^{43}V(0) \rightarrow V(0)
\end{align*}
This Massey product therefore detects the Toda bracket
\begin{equation}\label{eq:v_2^3b_1b_0^3}
\langle \beta_2 \beta_1^3\alpha_1[5], v_1, \td{\beta_3} \rangle.
\end{equation}
We
saw in 
Lemma~\ref{lem:dv_3h_0b_0} that $b_0^3b_1 = 0$ in $E_5(V(1))$.
Therefore, the Toda bracket is detected in a higher Adams filtration 
modulo an indeterminacy 
contained in the subgroup
$$ [\Sigma^{69}V(0), V(1)] \circ \td{\beta_3}. $$
We computed $\pi_{69}(V(1))$ and $\pi_{70}(V(1))$ using the AHSS
(\ref{AHSS:68-73}).  Both of these groups 
are rank $1$, generated by $x_{68}[1]$ and
$\beta_2^2\beta_1\alpha_1[5]$, respectively.
So we have $[\Sigma^{69}V(0), V(1)]$ is of rank $2$, generated 
by elements $\beta_2^2
\beta_1 \alpha_1[5]_1$ and $x_{68}[1]_0$.  
Here, the subscripts $0,1$ are used to indicate what cell of the
\emph{source} $V(0)$ the element is born on.
We will conclude that the
indeterminacy is zero using the $V(0)$-module structure of $V(1)$.
\begin{gather}
\label{eq:beta1_alpha1_indeterminacy}
\beta_2^2\beta_1\alpha_1[5]_1 \circ \td{\beta_3} = 
\beta_3\cdot\beta_2^2\beta_1\alpha_1[5]_0 = 0 \\
\label{eq:x68_beta3_indeterminacy}
x_{68}[1]_0 \circ \td{\beta_3} = (\beta_3[1] \cdot x_{68}[1])_0 = 0
\end{gather}
Equation~\ref{eq:beta1_alpha1_indeterminacy} 
follows from  the relation $\beta_3\beta_1 = 0$ in $\pi_*(S)$.

In Equation~\ref{eq:x68_beta3_indeterminacy}, 
we mean the product of $\beta_3[1] \in \pi_{43}(V(0))$ and $x_{68}[1] 
\in \pi_{69}(V(1))$ under the module map $\mu:V(0)\wedge V(1) \rightarrow
V(1)$.  The following diagram commutes.
$$
\xymatrix{
\pi_{43}(V(0))\otimes \pi_{69}(V(0)) \ar[r]_-{\mu} \ar[d] & \pi_{112}(V(0))
\ar[d]
\\
\pi_{43}(V(0)) \otimes \pi_{69}(V(1)) \ar[r]_-{\mu} & \pi_{112}(V(1))
}
$$
Since $x_{68}[1] \in \pi_{69}(V(1))$
is born on the $1$-cell, we may lift it to an
element $x_{68}[1] \in \pi_{69}(V(0))$.  Thus in order to show there is no
indeterminacy it suffices to show that $\beta_3[1] \cdot x_{68}[1] = 0 \in
\pi_{112}(V(0))$.
At this point we remind the reader that by our descriptions of the elements
$\beta_3[1]$ and $x_{68}[1]$ as elements in $\pi_*(V(1))$, we mean that
their images under the projection on to the top cell are $\beta_3$ and
$x_{68}$, respectively.  These elements are not necessarily uniquely
determined in $\pi_*(V(0))$, but any two representatives will differ by the
image of an element of $\pi_*(S)$ under the inclusion of the bottom cell of
$V(0)$.  Now $\pi_{69}(S) = 0$, so $x_{68}[1]$ is uniquely determined.
However, $\pi_{43}(S)$ is generated by $\alpha_{11}$, so any two elements of
$\pi_{43}(V(0))$ which project to $\beta_3$ on the top cell must differ by
$\pm \alpha_{11}[0] \in \pi_{43}(V(0))$.  
Under the product map $V(0)\wedge V(0) \rightarrow V(0)$, 
all homotopy carried by the
smash product of the $1$-cells is annihilated.  We will necessarily have
$\beta_3[1]\cdot x_{68}[1] = 0 \in \pi_{112}(V(0))$ for all possible
representatives of $\beta_3[1]$ if we can show
$$ \alpha_{11}[0] \cdot x_{68}[1] = 0 \in \pi_{112}(V(0)). $$
This is straightforward: in $\pi_{43}(V(0))$ we have $\alpha_{11}[0] = v_1^{10}
\alpha_1[0]$, and $\alpha_1 \cdot x_{68} = 0 \in \pi_*(S)$.
Thus the Toda bracket~\ref{eq:v_2^3b_1b_0^3} has no indeterminacy.

We conclude that the
Toda bracket \ref{eq:v_2^3b_1b_0^3} must be zero in $\pi_*(V(1))$
modulo higher Adams filtration.  
There
are only three elements in the correct range to kill the corresponding 
Massey product~\ref{eq:MasseyProduct}, and
these elements are $v_2^3v_3h_0b_0$, $\gamma_2 v_2^2$, and $v_3k_0h_0v_2^2$.  
In Lemma~\ref{lem:target9}, we proved that $v_2^3v_3h_0b_0$ supported a
non-trivial $d_2$.  The element $\gamma_2 v_2^2$ must support a non-trivial
$d_3$; this follows from the $S$-module pairing and the $d_3$ supported by
$v_2^2$ that was proved in Lemma~\ref{lem:dv_2^i}.  
Thus we must have
$$ d_5(v_3k_0h_0v_2^2) \doteq v_2^3b_1b_0^3 $$
and the lemma is proven after multiplying both sides by $b_0^3$.
\end{proof}

\begin{lem}\label{lem:target8}
In $E_2(V(1))$ there is a differential
$$ d_2(\eta_1 g_0 v_2^3 b_0^2) = \eta_1 g_0 k_0 b_0^3 h_1 \doteq k_1 b_0^6.$$
\end{lem}

\begin{proof}
This lemma follows immediately from the $S$-module pairing of Adams
spectral sequences. We know $b_0^2 \eta_1$ 
survives to $E_5(S)$ and we have computed
$d_2(g_0 v_2^3)$ in $E_2(V(1))$ as part of Lemma~\ref{lem:dg_0v_2^i}.
\end{proof}

We have established that every possible target of a differential supported
by $v_2^9$ is either the target of a shorter differential
or the source of a differential.  We may
conclude that $v_2^9$ is a permanent cycle.

\section{Proof that $v_2^9$ extends over $V(1)$}\label{sec:mainthmpf2}

In this section we will prove that if $v_2^9:S^{144} \rightarrow V(1)$ is a
map detected by the element $v_2^9 \in E_2(V(1))$, then it extends over
$V(1)$ to a self-map
$$ \td{v_2^9}: \Sigma^{144} V(1) \rightarrow V(1). $$
Applying Spanier-Whitehead duality, this is equivalent to finding an
extension corresponding to the dotted arrow in the diagram below.
$$\xymatrix{
& D(V(1)) \wedge V(1) \ar[d] \\
S^{144} \ar@{.>}[ur] \ar[r]_{v_2^9} & V(1) 
}$$
The  
inclusion of the wedge summand $Y_1$ of $D(V(1)) \wedge V(1)$ 
(Proposition~\ref{prop:splitting}) has the property that the composite
$$ Y_1 \hookrightarrow D(V(1)) \wedge V(1) \rightarrow V(1) $$ 
is just projection onto the top $V(1)$.
It therefore suffices to extend $v_2^9$ over the complex $Y_1$ 
as displayed below.
$$\xymatrix{
& Y_1 \ar[d] \\
S^{144} \ar@{.>}[ur] \ar[r]_{v_2^9} & V(1) \ar[d]^{\delta} \\
& \Sigma^{-5} V(1)
}$$
The vertical column forms a cofiber sequence where $\delta$ is given as the
composite
$$ V(1) \xrightarrow{\nu} \Sigma^{5}V(0) \xrightarrow{\beta_1}
\Sigma^{-5}V(0) \xrightarrow{\iota} \Sigma^{-5} V(1).$$
Here, $\nu$ is projection onto the top $V(0)$ and $\iota$ is inclusion of
the bottom $V(0)$.
Thus there is a solution to the extension problem if and only if the
composite $\delta \circ v_2^9 = 0$.  The map $\delta$ is also given by
the composite
$$\xymatrix{
V(1) \ar[r]^{\nu} \ar@/_2pc/[rr]_{\delta'} 
& \Sigma^{5}V(0) \ar[r]^{\iota} & \Sigma^{5} V(1) \ar[r]^{\beta_1} &
\Sigma^{-5} V(1)}$$
Here $\delta'$ is the geometric $v_1$-Bockstein.  In this section we will
prove
\begin{equation}\label{eq:delta'}
\beta_1 \cdot \delta'(v_2^9) = 0 
\end{equation}
from which it follows that $\delta \circ v_2^9 = 0$, and thus $v_2^9$
extends over $V(1)$.

Since $v_2^9$ has Adams filtration $9$, $\delta'(v_2^9)$ has Adams
filtration $\ge 9$.  We may calculate $\delta'$ in $\ext$ using the
periodic lambda algebra.  We have
$$ d(v_2^9) = v_1^9 \lambda_{27} $$
thus we have
$$ \xymatrix@R-2em{
\br{\Lambda}_{(1)} \ar[r]_{\nu_*} & \br{\Lambda}_{(0)} \ar[r]_{\iota_*} &
\br{\Lambda}_{(1)} \\
v_2^9 \ar@{|->}[r] & v_1^8 \lambda_{27} \ar@{|->}[r] & 0
} $$
and so $\delta'_*(v_2^9) = 0$ viewed as an element of $E_2(V(1))$.  We may
conclude that the Adams filtration of $\delta'(v_2^9)$ is $\ge 10$.  Our
strategy to prove that $\beta_1 \cdot \delta'(v_2^9) = 0$ 
is to make a list of all of the
elements in $E_2(V(1))$ in Adams filtration greater than or equal to $10$ in
the $139$-stem.  We then will prove that each of these elements is either
not a permanent cycle, or is killed by a differential or at least has the
property that composition with $\beta_1$ is zero.  A list of the
elements, as well as the lemmas that deal with it, is given below.
\begin{center}
\begin{tabular}{c|c|c}
\hline
AF & Element & Lemma \\
\hline \hline
26 & $h_0v_2b_0^{12}$ & \ref{lem:target2} \\
\hline
25 & $k_0h_0b_0^{11}$ & \ref{lem:target2} \\
\hline
20 & $g_0h_0v_2^3b_0^7$ & \ref{lem:target10} \\
\hline
19 & $b_1h_0v_2^2b_0^7$ & \ref{lem:target11} \\
\hline
15 & $v_2^6h_0b_0^4$ & \ref{lem:stem139}, \ref{lem:target12} \\
 & $\eta_1v_2^2b_0^5$ & \ref{lem:stem139}, \ref{lem:target13} \\
 & $v_3h_1v_2b_0^6$ & \ref{lem:stem139}, \ref{lem:target13} \\
\hline
14 & $\eta_1k_0v_2b_0^4$ & \ref{lem:stem139}, \ref{lem:target14} \\
 & $v_2^5k_0h_0b_0^3$ & \ref{lem:stem139}, \ref{lem:target15} \\
 & $v_3k_0h_1 b_0^5$ & \ref{lem:stem139}, \ref{lem:target14}, 
 \ref{lem:target15} \\
\hline
13 & $\langle k_1, h_0, h_0 \rangle b_0^5$ & \ref{lem:target16} \\
\hline
10 & $v_3h_0b_0^2v_2^4$ & \ref{lem:target17} \\
\hline
\end{tabular}
\end{center}

\begin{lem}\label{lem:target10}
If $x$ is an element of $E_r(V(1))$ whose Hurewicz image 
$h(x) \in E_r(eo_2 \wedge
V(1))$ supports a non-trivial differential, then $x$ cannot be a permanent
cycle.
\end{lem}

\begin{proof}
This is obvious; $h$ is a map of spectral sequences.
\end{proof}

We shall need the following lemma.

\begin{lem}\label{lem:stem68}
There is the following pattern of differentials in the ASS for
$\pi_*(V(1))$ in the vicinity of the $68$-stem.
$$ \xymatrix@R-2em{
\underline{\text{Stem $67$}} & \underline{\text{Stem $68$}} 
& \underline{\text{Stem $69$}} \\
v_2b_0^4h_1 & b_0^5g_0 & b_0^5v_2h_0 \\
b_1b_0^3h_0 & b_0^2v_2^3 \ar[l] & b_0^4h_0k_0 \\
\eta_1b_0 & b_0v_2^2k_0 \ar@{.>}[uul] & b_0h_1v_2^3 \ar[uul] \\
v_2^4h_0 & v_2 k_0^2 & h_0b_1v_2^2 \\
& \eta_1 h_1 \ar@{.>}[uuuul] & g_1h_1}
$$
Here only one of the dotted differentials occurs.
\end{lem}

\begin{rmk}
We  will find later
(see the proof of Lemma~\ref{lem:stem139}) that in fact $\eta_1h_1$ must be
a permanent cycle, and thus the dotted differential supported on
$b_0v_2^2k_0$ must be the non-trivial one.
\end{rmk}

\begin{proof}
We will first deduce the following portion of the ASS chart for
$\pi_*(V(1))$.
$$\xymatrix@R-2em{
\underline{\text{Stem $56$}} & \underline{\text{Stem $57$}} & 
\underline{\text{Stem $58$}} & 
\underline{\text{Stem $59$}} \\
b_0^4v_2 & b_0^3h_1v_2 & b_0^4g_0 & v_2b_0^4h_0 \\
b_0^3k_0 & b_1b_0^2h_0 & b_0v_2^3 \ar[l] & k_0b_0^3h_0 \\
& \eta_1 \ar[ul] & v_2^2k_0 \ar@{.>}[uul] & h_1v_2^3 \ar@{.>}[uul]^{(2)} \\
& & g_1 \ar@{.>}[uuul]^{(1)} 
}$$
The differential on $\eta_1$ follows from the differential in $E_5(S)$.
The differential on $b_0v_2^3$ is a consequence of Lemma~\ref{lem:dv_2^i}.  
Since the AHSS calculation (\ref{AHSS:55-58})
told us that $\pi_{57}(V(1)) = 0$, we may
conclude that one of the dotted differentials $(1)$ exists.  Also, we have
shown that $\pi_{58}(V(1))$ has rank $1$, hence something must kill
$b_0^4g_0$.  Both $v_2b_0^4h_0$ and $k_0b_0^3h_0$ are permanent cycles, so
the only candidate to support the dotted differential $(2)$ is $h_1v_2^3$.
Note that this differential is present in the ASS for $eo_2\wedge V(1)$.

Moving up to the vicinity of the $68$-stem of the ASS,
the two solid differentials in the statement of the lemma 
are propagated by $b_0$ multiplication.
Our AHSS calculations (\ref{AHSS:63-68}) 
tell us that $\pi_{67}(V(1))=0$.  Now 
$v_2k_0^2$ must be a permanent cycle since  
$v_2k_0$ is a permanent
cycle for dimensional reasons.  Since $v_2b_0^4h_1$ must vanish, one of
the dotted differentials must occur.  
We have computed $\pi_{68}(V(1))$ to
be of rank $2$, so there can be no more differentials originating from the
$69$ stem. 
\end{proof}

\begin{lem}\label{lem:target11}
The element $b_1h_0v_2^2b_0^7 \in E_6(V(1))$ must be zero.
\end{lem}

\begin{proof}
In Lemma~\ref{lem:stem68} we showed that $v_2^2 b_1 h_0 \in E_r(V(1))$ 
is a permanent cycle.  The
result then follows from the fact that in $E_7(S)$ there is a relation
$b_0^6 = 0$.  
\end{proof}

We are now in a position to prove the differential promised in
Remark~\ref{rmk:v_2^5}.  We will need this differential later in this
paper.

\begin{lem}\label{lem:v_2^5}
In $E_2(V(1))$, there is a non-trivial differential
$$d_2(v_3g_0) \doteq k_0h_1v_2^2 = b_1h_0v_2^2.$$
\end{lem}

\begin{proof}
Our AHSS computations (\ref{AHSS:68-73})
show that $\pi_{69}(V(1))$ has rank $1$.
In the proof of Lemma~\ref{lem:stem68}, we computed all of the Adams
differentials in $E_*(V(1))$ supported in the $69$ stem.  That data is used
to compute the following portion of the ASS for $\pi_*(V(1))$.
$$\xymatrix@R-2em{
\underline{\text{Stem $68$}} & \underline{\text{Stem $69$}} 
& \underline{\text{Stem $70$}} \\
b_0^5g_0 & b_0^5v_2h_0 & b_0^7 \\
b_0^2v_2^3 & b_0^4h_0k_0 & b_0^2g_0v_2^2 \ar[l] \\
b_0v_2^2k_0 & b_0h_1v_2^3 \ar[uul] & b_0^2v_2b_1 \ar[uul] \\
v_2k_0^2 & h_0b_1v_2^2 & h_0b_0\eta_1 \\
\eta_1h_1 & g_1h_1 & v_3g_0 \ar@{.>}[ul]
}$$
The differential supported by $b_0^2g_0v_2^2$ is a consequence of
Lemma~\ref{lem:dg_0v_2^i}.  The differential supported by $b_0^2v_2b_1$ is
a consequence of the Toda differential on $b_1$ in $E_5(S)$ and the
$S$-module structure of $V(1)$.  One more element in the $69$ stem must be
the target of an Adams differential, and the only possibility is the dotted
differential.  This is the differential we wanted.
\end{proof}

There is an elaborate pattern of activity in Adams filtration $13$-$15$ in
the $139$-stem.  
We would like to understand which elements in this range
of filtrations are permanent cycles, and which aren't, and this is
accomplished in Lemma~\ref{lem:stem139}.  As a consequence,
certain linear combinations of the possible obstructions to the extension
of $v_2^9$ will be eliminated.  First we need the following differential.

\begin{lem}\label{lem:dv_3h_1v_2}
In $E_3(V(1))$, there is a non-trivial differential
$$ d_3(v_3h_1v_2) \doteq b_0\eta_1h_1. $$
\end{lem}

\begin{proof}
In Lemma~\ref{lem:target4}, we showed that $v_3h_1$ supports a non-trivial
$d_4$ and thus is a $d_3$ cycle. We now apply the product formula
(\ref{ProductRule}) to deduce the differential on $v_3h_1\cdot v_2$.
Computer assisted lambda algebra calculations yield the formula
$$ \br{v_3h_1} = -g_1 \in E_2(V(0)). $$
We wish to compute $\brr{v_3h_1}$, which is obtained by computing the
Bockstein on $g_1$.  Now $g_1$ is a $d_1$-cycle when it is considered as an
element of $E_1(S^0)$, so we need to compute the Adams differential of
$g_1$ in $E_*(S^0)$ to get this Bockstein.  We use the main theorem of Bruner
(\cite[VI.1.1]{Hinfty}) to understand the relationship between Steenrod
operations in $\ext$ and differentials in the ASS.  Using Bruner's formula,
we may compute
\begin{equation}\label{eq:betaP^0}
d_2(g_1) = d_2(P^0(g_0)) = v_0\cdot \beta P^0(g_0) \doteq v_0\eta_1.
\end{equation}
Here we should point out that our indexing of the Steenrod operations is
different from that of Bruner's, but conforms to the perhaps more common
indexing as given in \cite{May}.  The element $v_0$ detects the degree $p$
map on spheres.  
The Steenrod operation $\beta P^0(g_0) \doteq \eta_1$ is computed using the
May spectral sequence for $H^*(P_*)$.  Specifically,  
The element $g_0$ is detected by by $\pm h_{2,0}h_{1,0}$.  On the May
$E_2$-term, we compute (using the Cartan formula)
$$ \beta P^0(h_{2,0}h_{1,0}) = \beta P^0(h_{2,0})\cdot P^0(h_{1,0}) -
P^0(h_{2,0}) \cdot \beta P^0(h_{1,0}) = b_{2,0}h_{1,1} - h_{2,1}b_{1,0} $$
and $\pm(b_{2,0}h_{1,1} - h_{2,1}b_{1,0})$ detects $\eta_1$.

We conclude from Equation~\ref{eq:betaP^0}
that $\brr{v_3h_1} \doteq \eta_1$.  We then apply the product
formula, keeping in mind that $\br{v_2} = h_1$ and $\brr{v_2} = b_0$, 
and get
$$ d_3(v_3h_1 \cdot v_2) = -b_0 (\pm \eta_1h_1 - g_1b_0) \doteq
b_0\eta_1h_1. $$
\end{proof}

\begin{lem}\label{lem:stem139}
We have the following differentials in the $139$ stem of the ASS for
$V(1)$ in Adams filtrations $13$--$15$.  
\begin{align*}
d_3(\eta_1v_2^2b_0^5) & \doteq b_0^7h_1\eta_1 \\
d_3(v_3h_1b_0^6v_2) & \doteq b_0^7h_1\eta_1 \\
d_3(v_2^6h_0b_0^4) & \doteq b_0^7h_1\eta_1 \quad \text{or $v_2^6h_0$ is a
permanent cycle} \\
d_5(v_2^5k_0h_0b_0^3) & \doteq b_0^7k_0^2v_2 \\
d_5(\eta_1k_0v_2b_0^4) & \doteq b_0^7k_0^2v_2 \\
d_5(v_3k_0b_0^5h_1) & \doteq b_0^7k_0^2v_2 
\end{align*}
An $\FF_3$-basis of permanent cycles in this range
of Adams filtration is listed below.
\begin{gather*}
v_2^6h_0b_0^4 + a_1 (\eta_1v_2^2b_0^5) + a_2 (v_3h_1b_0^6v_2) \\
\eta_1v_2^2b_0^5 \pm v_3h_1b_0^6v_2 \\
v_2^5k_0h_0b_0^3 \pm v_3k_0b_0^5h_1 \\
\eta_1k_0v_2b_0^4 \pm v_3k_0b_0^5h_1 \\
\langle k_1, h_0, h_0 \rangle b_0^5
\end{gather*}
Here $a_1$ and $a_2$ are elements of $\FF_3$.
We are unable to determine the $\pm$ signs or the coefficients $a_i$.  
A diagram of this
portion of the ASS chart is given below for the reader's convenience.
$$\xymatrix@R-2.5em@C+2em{
\underline{\text{AF}} & \underline{\text{Stem $138$}} &
\underline{\text{Stem $139$}} & \\
19 & k_0^2v_2b_0^7 && \\
\ar@{.}[rrr] &&& \\
18 & b_0^7 h_1 \eta_1 &&& \\
\ar@{.}[rrr] &&& \\
15 && v_2^6h_0b_0^4 \ar@{-->}[uul] & \\
   && \eta_1v_2^2b_0^5 \ar@<+.5ex>[uuul] & \\
   && v_3h_1b_0^6v_2 \ar@<1ex>[uuuul] & \\
\ar@{.}[rrr] &&& \\
14 && v_2^5k_0h_0b_0^3 \ar@<-1ex> `l[u] `[uuuuuuuul] [uuuuuuuul] & \\
   && \eta_1k_0v_2b_0^4 \ar `l[u] `[uuuuuuuuul] [uuuuuuuuul] & \\
   && v_3k_0b_0^5h_1 \ar@<+1ex> `l[u] `[uuuuuuuuuul] [uuuuuuuuuul] & \\
\ar@{.}[rrr] &&& \\
13 && \langle k_1, h_0, h_0 \rangle b_0^5 
}$$
\end{lem}

\begin{proof}
The method of proof is to divide these elements by a maximal power of
$b_0$, and then multiply by $b_0$ successively until all of the elements in
question are present in the $E_2$ term.  We begin with the vicinity of the
$79$ stem.  
In our AHSS calculations (\ref{AHSS:77-80}), we computed $\pi_{78}(V(1))$
and $\pi_{79}(V(1))$ modulo one differential which we were unable to
determine (this is the dotted differential 
labeled (1) in (\ref{AHSS:77-80})).
We conclude that either $\pi_{78}(V(1))$ has rank $2$ and $\pi_{79}(V(1))$
has rank $1$, or $\pi_{78}(V(1))$ has rank $1$ and $\pi_{79}(V(1))$ is
trivial.
We will see shortly however that the latter is the case, i.e. that the dotted
differential $(1)$ must exist.
We display below the ASS in the same range.
$$\xymatrix@R-2em{
\underline{\text{Stem $77$}} & \underline{\text{Stem $78$}} 
& \underline{\text{Stem $79$}} & \underline{\text{Stem $80$}} \\
v_2b_0^5h_1 & b_0^6g_0 & v_2b_0^6h_0 & b_0^8 \\
b_1b_0^4h_0 & v_2^3b_0^3 \ar[l] & k_0b_0^5h_0 & b_0^3g_0v_2^2 \ar[l] \\
b_0^2\eta_1 & v_2^2k_0b_0^2 \ar[uul]& v_2^3b_0^2h_1 \ar[uul] & b_0^3v_2b_1 
\ar[uul] \\
v_2^4h_0b_0 & b_0v_2k_0^2 & b_1v_2^2b_0h_0 & \eta_1b_0^2h_0 \\
v_2^3k_0h_0 & \eta_1h_1b_0 & v_3h_1v_2 \ar[l] & b_0v_3g_0 \ar[ul] \\
& v_3k_0 \ar@{.>}[ul]^{(2)} & & v_2^5 \ar[uul]
}$$
By comparing with the ASS chart in the vicinity of the $68$-stem in the
proof of Lemma~\ref{lem:stem68}, we see that
the elements $b_0v_2k_0^2$, $v_2b_0^6h_0$, $k_0b_0^5h_0$, and
$b_1v_2^2b_0h_0$ are permanent cycles, and the elements 
$v_2^3b_0^3$ and $v_2^3b_0^2h_1$ support the
indicated differentials.  The differential on $v_3h_1v_2$ was the subject
of Lemma~\ref{lem:dv_3h_1v_2}.  In Lemma~\ref{lem:stem68}, we were unable
to determine whether $v_2b_0^4h_1$ was killed by $b_0v_2^2k_0$ or
$\eta_1h_1$.  Since $\eta_1h_1b_0$ is the target of a differential, this
ambiguity is now resolved: $d_4(v_2^2k_0b_0^2) \doteq v_2b_0^5h_1$.
The differential supported by $b_0^3g_0v_2^2$ was established in
Lemma~\ref{lem:dg_0v_2^i}.  The differential supported by $b_0^3v_2b_1$
follows from the Toda differential on $b_1$ in $E_*(S^0)$.  The
differentials supported by $b_0v_3g_0$ and $v_2^5$ were proven in
Lemmas~\ref{lem:v_2^5} and \ref{lem:dv_2^i}, respectively.  There is
nothing remaining in stem $79$, so we conclude that $\pi_{79}(V(1))$ is 
trivial.  Therefore the dotted differential $(1)$ exists in the 
AHSS chart (\ref{AHSS:77-80}).  The AHSS chart now tells us that 
$\pi_{78}(V(1))$ is of rank
$1$, and the only way for the ASS to produce the same answer is for there
to exist the dotted differential $(2)$ since neither $b_0^2\eta_1$ nor
$v_2^4h_0b_0$ are permanent cycles.

We now multiply everything by $b_0$ and move into the vicinity of the $89$
stem.  
The ASS chart is displayed below.
$$\xymatrix@R-2em{
\underline{\text{Stem $87$}} & \underline{\text{Stem $88$}} 
& \underline{\text{Stem $89$}} & \underline{\text{Stem $90$}} \\
b_0^6h_1v_2 & b_0^7g_0 & v_2b_0^7h_0 & b_0^9 \\
b_1b_0^5h_0 & v_2^3b_0^4 \ar[l] & k_0b_0^6h_0 & v_2^2b_0^4g_0 \ar[l] \\
\eta_1b_0^3 & v_2^2k_0b_0^3 \ar[uul] & v_2^3h_1b_0^3 \ar[uul] & b_1v_2b_0^4
\ar[uul] \\
v_2^4b_0^2h_0 & b_0^2v_2k_0^2 & b_1v_2^2b_0^2h_0 & \eta_1b_0^3h_0 \\
b_0h_0v_2^3k_0 & b_0^2h_1\eta_1 & v_2^2\eta_1 \ar[l] & b_0^2v_3g_0 \ar[ul] \\
v_3h_2 & v_3k_0b_0 \ar[ul] & b_0v_2v_3h_1 \ar[ul]  & v_2^5b_0 \ar[uul] \\
&& h_1v_3k_0 \ar[uuul] & v_2^4k_0 \\
&& \langle k_1, h_0, h_0 \rangle & v_2^2g_1 \\
&&& v_3h_2h_0
} $$
All of the indicated differentials follow from our computations near the
$79$ stem except for the one supported by $h_1v_3k_0$.  For that we
consider the image under projection onto the top cell of $V(1)$. We saw in
the proof of Lemma~\ref{lem:dv_3h_1v_2} that $\brr{v_3h_1} \doteq \eta_1$.  
We have
$$ \brr{d(v_3h_1k_0)} = d(\brr{v_3h_1k_0}) \doteq d(\eta_1k_0) \doteq
k_0^2b_0^3 = \brr{b_0^2v_2k_0^2}. $$
This can only happen if $d_5(v_3h_1k_0) = k_0^2b_0^2v_2$.  Our AHSS
calculations (\ref{AHSS:87-90})
tell us that $\pi_{89}(V(1))$ is of rank $2$, so there can be
no more differentials originating from the $90$ stem.

We now multiply once more by $b_0$ and arrive in the crucial region around
the $99$ stem.  
Our AHSS computations (\ref{AHSS:97-100}) tell us that the rank of 
$\pi_{98}(V(1))$ is $3$ and the rank of
$\pi_{99}(V(1))$ is $4$.  We turn now to the ASS.
$$\xymatrix@R-2em{
\underline{\text{Stem $97$}} & \underline{\text{Stem $98$}} 
& \underline{\text{Stem $99$}} & \underline{\text{Stem $100$}} \\
b_0^7v_2h_1 & b_0^8g_0 & v_2b_0^8h_0 & b_0^{10} \\
b_1b_0^6h_0 & v_2^3b_0^5 \ar[l] & k_0b_0^7h_0 & b_0^5g_0v_2^2 \ar[l] \\
b_0^4\eta_1 & b_0^4v_2^2k_0 \ar[uul] & v_2^3b_0^4h_1 \ar[uul] &
b_0^5v_2b_1 \ar[uul] \\
v_2^4b_0^3h_0 & b_0^3k_0^2v_2 & v_2^2b_1b_0^3h_0 & b_0^4h_0\eta_1 \\
v_2^3b_0^2h_0k_0 & b_0^3h_1\eta_1 & v_2v_3b_0^2h_1 \ar[l] 
& b_0^3v_3g_0 \ar[ul] \\
v_2^2v_3h_0b_0 & b_0^2v_3k_0 \ar[ul] 
& h_0v_2^6 \ar@{-->}[ul]|-{(4)}
& v_2^5b_0^2 \ar[uul] \\
v_2\gamma_2 & g_0v_2^5 \ar@{-->}[uul]|-{(3)} & b_0v_2^2\eta_1 \ar[uul] &
v_2^4k_0b_0 \\
b_0v_3h_2 & v_2^4b_1 & v_2\eta_1k_0 \ar@<-1ex> `l[u] `[uuuul] [uuuul] 
& v_2^3k_0^2 \\
& v_2k_1 & b_0h_1v_3k_0 \ar@<-.5ex> `l[u] `[uuuuul] [uuuuul] 
& v_2^2h_1\eta_1 \\
& & \langle x_{92}, h_0, h_0 \rangle & h_0v_2\gamma_2
}$$
Aside from the differentials supported by (possibly) $g_0v_2^5$,
$h_0v_2^6$, and $v_2\eta_1k_0$, all of the differentials displayed follow
from our calculations near the $89$ stem.  If we had only the
differentials arising from $b_0$ multiplication 
on elements in the vicinity of the $89$
stem, we would have created groups of the 
correct rank in the $98$ and $99$ stem
as predicted from the AHSS.  Therefore any additional differentials must
preserve the rank of the $E_{\infty}$ term.  We easily see that
$d_5(v_2\eta_1k_0) = b_0^3k_0^2v_2$ from the differential on $\eta_1$ and
the $S$-module structure of $V(1)$.

The problem is that $h_0v_2^6$
could support a $d_4$ killing $b_0^3k_0^2v_2$, and this would make both
$v_2\eta_1k_0$ and $b_0h_1v_3k_0$ into permanent cycles.  We claim that 
this cannot happen.  For suppose that $d_4(h_0v_2^6) \doteq b_0^3k_0^2v_2$.
Then there is no linear combination of elements containing $h_0v_2^6$ which
is a permanent cycle.
Our AHSS calculation indicates there is some element in $\pi_{99}(V(1))$
such that its image in $\pi_{93}(S)$ under projection onto the top cell of
$V(1)$ is $\beta_6\alpha_1$.  
The only element which can account for this is $h_0v_2^6$.  Therefore, if
$h_0v_2^6$ is not a permanent cycle, it must support a $d_3$ killing
$b_0^3h_1\eta_1$.  This possibility is indicated by the dashed differential
$(4)$.

Also, we cannot eliminate the possibility that $d_3(g_0v_2^5) \doteq
v_2^3b_0^2h_0k_0$, since $g_0v_2^5$ is in the same Adams filtration as
$b_0^2v_3k_0$.  
This possibility is indicated with the dashed differential $(3)$.
Upon taking the pattern of differentials supported by the
$99$ stem and multiplying by $b_0^4$, we get the promised pattern of
differentials supported in the $139$ stem.
\end{proof}

\begin{lem}\label{lem:target12}
Choose the correct coefficients $a_i \in \FF_3$ so that
$v_2^6h_0b_0^4+a_1(\eta_1v_2^2b_0^5)+a_2(v_3h_1b_0^6v_2)$ is a permanent
cycle.  The composite of any element that this permanent cycle detects
with $\beta_1$ must be null.
\end{lem}

\begin{proof}
The element $v_2^6h_0+a_1(\eta_1v_2^2b_0)+a_2(v_3h_1b_0^2v_2)$ is a
permanent cycle by Lemma~\ref{lem:stem139}.  Now apply
Corollary~\ref{cor:beta_1_order}.  Comparing with the elements present 
in the $149$
stem of $E_2(V(1))$ we see that there is no possibility of a hidden
$\beta_1$ extension.
\end{proof}

\begin{lem}\label{lem:target13}
Choose the correct sign so that $\eta_1v_2^2b_0^5 \pm v_3h_1b_0^6v_2$ is a
permanent cycle.
This element must be the target of a differential.
\end{lem}

\begin{proof}
The element $\eta_1v_2^2\pm v_3h_1b_0v_2$ is a
permanent cycle (Lemma~\ref{lem:stem139}).  Now apply
Corollary~\ref{cor:beta_1_order}.
\end{proof}

\begin{lem}\label{lem:target14}
Choose the correct sign so that $\eta_1k_0v_2b_0^4\pm v_3k_0b_0^5h_1$ is a
permanent cycle.  The composite of any element that this permanent cycle 
detects with $\beta_1$ must be null.
\end{lem}

\begin{proof}
Again, apply Corollary~\ref{cor:beta_1_order}.  Comparing with the elements 
present 
in the $149$
stem of $E_2(V(1))$ we see that there is no possibility of a hidden
$\beta_1$ extension.
\end{proof}

\begin{lem}\label{lem:target15}
Choose the correct sign so that $v_2^5k_0h_0b_0^3\pm v_3k_0b_0^5h_1$ is a
permanent cycle.  This element must be the target of a differential. 
\end{lem}

\begin{proof}
In Lemma~\ref{lem:v_2^5}, it was established that $v_2^5 \pm v_3g_0b_0$ was
a permanent cycle, therefore the element 
$$ h_0 \cdot (v_2^5 \pm v_3g_0b_0) = h_0v_2^5\pm v_3h_1b_0^2 $$
is a permanent cycle.  Now $d_5(\eta_1) \doteq k_0b_0^3$ in the ASS for
$\pi_*(S)$.  Using the $S$-module structure of $V(1)$, we deduce that there
is a differential
$$ d_5(\eta_1\cdot(h_0v_2^5\pm v_3h_1b_0^2)) \doteq 
v_2^5k_0h_0b_0^3\pm v_3k_0b_0^5h_1. $$
Note that $v_2^5k_0h_0b_0^3\pm v_3k_0b_0^5h_1$ might be the target of a
shorter differential, in which case the conclusion of the lemma is still
satisfied.
\end{proof}

\begin{lem}\label{lem:target16}
The element $\langle k_1, h_0,
h_0 \rangle b_0^5$ must be the target of a differential.
\end{lem}

\begin{proof}
This follows immediately from Corollary~\ref{cor:beta_1_order} and the fact
that $\langle k_1, h_0, h_0 \rangle$ is a permanent cycle
(Lemma~\ref{lem:stem139}).
\end{proof}

\begin{lem}\label{lem:target17}
In $E_2(V(1))$ there is a non-trivial differential
$$ d_2(v_3h_0b_0^2v_2^4) \doteq v_2^2h_0\eta_1k_0b_0^2. $$
\end{lem}

\begin{proof}
In Lemma~\ref{lem:dv_3h_0b_0}, we showed that $d_2(v_3h_0b_0) = 0$.  In
Lemma~\ref{lem:dv_2^i}, we showed that 
$$d_2(v_2^4) = -b_0k_0h_1v_2 = -b_0b_1h_0v_2. $$
Computer assisted lambda algebra computations reveal that in $E_2$, we have
\begin{align*}
\br{v_3h_0b_0} & = -h_0\eta_1 \\
\brr{v_3h_0b_0} & = 0
\end{align*}
We also have $\br{v_2^4} = h_1v_2^3$ and $\brr{v_2^4} = 0$ in $E_2$.
We will now use the product rule (\ref{ProductRule}) to compute the
differential on $(v_3h_0b_0) \cdot v_2^4$.  
The term 
$$ G(v_3h_0b_0, v_2^4) = \br{v_3h_0b_0} \cdot \brr{v_2^4} -
\brr{v_3h_0b_0} \cdot \br{v_2^4} $$ 
is trivial in $E_2$ because $\brr{v_2^4}$ and $\brr{v_3h_0b_0}$ are trivial
in $E_2$.
The following manipulations are
made possible with the hidden extension $h_0\cdot(v_3h_0b_0) \doteq
b_1g_0v_2$.
\begin{align*}
d_2((v_3h_0b_0)\cdot v_2^4) & = (v_3h_0b_0)\cdot(b_0b_1h_0v_2) \\
& \doteq b_1^2g_0b_0v_2^2 \\
& \doteq v_2^2h_0\eta_1k_0b_0
\end{align*}
Now multiply by $b_0$.
\end{proof}


\begin{thebibliography}{99}

\bibitem{Adams}
J.F.~Adams, Stable Homotopy and Generalised Homology, University of Chicago
Press, Chicago, 1974.

\bibitem{Aubry}
M.~Aubry.  Calcul de groupes d'homotopie stables mod. $p$ par la suite 
spectrale
d'Adams-Novikov, C. R. Acad. Sc. Paris, S\'erie A 288 (1979) 587 -- 590.

\bibitem{Hinfty}
R.R.~Bruner, J.P.~May, M.~Steinberger, J.~McClure, $H_{\infty}$ Ring
Spectra and their Applications, Lecture Notes Math., Vol. 1176, Springer,
Berlin, 1980.

\bibitem{EKMM}
A.~Elmendorf, I.~Kriz, M.~Mandell, J.~P.~May, Rings, Modules, and Algebras
in Stable Homotopy Theory, Mathematical Surveys and Monographs, Vol. 47, 
Amer. Math. Soc., 1997.

\bibitem{GoerssHennMahowald}
P.~Goerss, H.-W.~Henn, M.~Mahowald, The homotopy of $L_2 V(1)$ for the
prime $3$, to appear in Skye Topology Conf. Proc.,  preprint
available from http://www.hopf.math.purdue.edu (2002).

\bibitem{GoerssHennMahowaldRezk}
P.~Goerss, H.-W.~Henn, M.~Mahowald, C.~Rezk, A resolution of the
$K(2)$-local sphere., preprint available from 
http://www.hopf.math.purdue.edu (2002).

\bibitem{Gray}
B. Gray, The periodic lambda algebra, Fields Institute Comm. 19 (1998),
93--101.

\bibitem{HopkinsMahowald}
M.J.~Hopkins, M.~Mahowald, From elliptic curves to homotopy theory,
preprint available from
http://www.hopf.math.purdue.edu (1998).

\bibitem{HopkinsMiller}
M.~Hopkins, H.~Miller, MIT Lecture Notes, Preprint.

\bibitem{HopkinsSmith}
M.J.~Hopkins, J.H.~Smith, Nilpotence and stable homotopy theory II., Ann.
of Math. (2) 148 (1998), 1--49. 

\bibitem{James}
I.M.~James, Reduced product spaces., Ann. of Math. (2) 62 (1955), 170--197.

\bibitem{Mahowald}
M.~Mahowald, Ring spectra which are Thom complexes, Duke Math. J. 46
(1979), 549--559.

\bibitem{May}
J.P.~May, A general algebraic approach to Steenrod operations, The Steenrod 
Algebra and its Applications (Proc. Conf. to Celebrate N.E.~Steenrod's 
Sixtieth Birthday, Battelle Memorial Inst., Columbus, Ohio, 1970),
Lecture Notes in Math., Vol. 168, Springer, Berlin, 1970, pp. 153--231.

\bibitem{Miller}
H.~Miller, On relations between Adams spectral sequences, with an
application to the stable homotopy of a Moore space, J. Pure and Applied
Alg. 20 (1981), 287--312.

\bibitem{Oka}
S.~Oka, Note on the $\beta$-family in stable homotopy of spheres at the
prime $3$, Kyushu Mem. Fac. Sci. Ser. A 35 (1981), 367--373.

\bibitem{Ravenel}
D.~Ravenel, Complex Cobordism and the Stable Homotopy Groups of Spheres,
Pure and Applied Math., Academic Press, 1986. 

\bibitem{Rezk}
C.~Rezk, Notes on the Hopkins-Miller theorem, Homotopy Theory via Algebraic
Geometry and Group Representations, Contemp. Math. 220 (1998), 313--366.

\bibitem{Shimomura:beta_1}
K.~Shimomura, On the action of $\beta_1$ in the stable homotopy of spheres
at the prime $3$, Hiroshima Math. J. 30 (2000), 345--362.

\bibitem{Shimomura:V(1)}
K.~Shimomura, The homotopy groups of the $L_2$-localized Toda-Smith
spectrum $V(1)$ at the prime $3$, Trans. Amer. Math. Soc. 349 (1997),
1821--1850.

\bibitem{Shimomura:V(0)}
K.~Shimomura, The homotopy groups of the $L_2$-localized mod $3$ Moore
spectrum, J. Math. Soc. Japan 52 (2000), 65--90.

\bibitem{Silverman}
J.~Silverman, The Arithmetic of Elliptic Curves, Graduate Texts in Math.,
Vol. 151, Springer-Verlag, New York, 1986.

\bibitem{Tangora}
M.~Tangora, Computing the Homology of the Lambda Algebra, Mem. of the
Amer. Math. Soc. 337, Amer. Math. Soc., Providence, 1985.

\bibitem{Thomas}
E.~Thomas, Whitney-Cartan product formulae, Math. Zeit. 118 (1970),
115--138.

\bibitem{Toda}
H.~Toda, An important relation in the homotopy groups of spheres, Proc.
Japan Acad. 43 (1968), 893--942.

\end{thebibliography}
\end{document}